%% file: main.tex
\documentclass[preprint,12pt]{elsarticle2}
\usepackage[english]{babel}
\usepackage[T1]{fontenc}
\usepackage{amsmath,amsthm}
\usepackage{amssymb}
\usepackage{amsfonts}
\usepackage{appendix}
\usepackage{geometry}
\usepackage{bbold}
\usepackage{mathrsfs}
\usepackage{graphicx}
\usepackage{tikz}
\usepackage{array}
\usepackage{amssymb}
\usepackage{natbib}

\journal{Mathematics and Computers in Simulation}
\usetikzlibrary{shapes.arrows,chains,patterns,matrix,arrows}
\newcommand*{\crochet}[1]{\left<#1\right>}

\newcommand*{\indi}[3]{#1 \mathbb{1}_{#3 > 0} + #2 \mathbb{1}_{#3 < 0}}

\newcommand*{\derivee}[2]{\dfrac{\mathrm{d}#1}{\mathrm{d}#2}}
\newcommand*{\sint}[4]{\int_{#1}^{#2} #3 \mathrm{d}#4}

\DeclareMathOperator{\sign}{sgn}

\DeclareMathOperator{\Dp}{\mathrm{F}_+}
\DeclareMathOperator{\Dm}{\mathrm{F}_{-}}

\DeclareMathOperator*{\argmin}{\arg\!\min}

\newtheorem{prop}{Proposition}[section]

\newtheorem{defi}{Definition}[section]
\newtheorem{rema}{Remark}[section]

\begin{document}
\begin{frontmatter}
\title{An asymptotic preserving kinetic scheme for the M1 model of linear transport}
\author[1]{Jean-Luc Feugeas}
\ead{jean-luc.feugeas@u-bordeaux.fr}

\author[2]{Julien Mathiaud}
\ead{julien.mathiaud@univ-rennes.fr}

\author[3]{Luc Mieussens}
\ead{Luc.Mieussens@u-bordeaux.fr}

\author[1]{Thomas Vigier\corref{cor1}}
\ead{thomas.vigier@u-bordeaux.fr}

\affiliation[1]{organization={CELIA, University of Bordeaux, CNRS, CEA},
            addressline={UMR 5107}, 
            city={Talence},
            postcode={F-33405},
            country={France}}
\affiliation[2]{organization={University of Rennes, CNRS, IRMAR},
            addressline={UMR 6625}, 
            city={Rennes},
            postcode={F-35000},
            country={France}}
\affiliation[3]{organization={University of Bordeaux, Bordeaux INP, CNRS, IMB},
            addressline={UMR 5251}, 
            city={Talence},
            postcode={F-33400},
            country={France}}
\cortext[cor1]{Corresponding author}

\begin{abstract}
Moment models with suitable closure can lead to accurate and computationally efficient solvers for particle transport. Hence, we propose a new asymptotic preserving scheme for the M1 model of linear transport that works uniformly for any Knudsen number. Our idea is to apply the M1 closure at the numerical level to an existing asymptotic preserving scheme for the corresponding kinetic equation, namely the Unified Gas Kinetic scheme (UGKS) originally proposed in \cite{kxu2010ugks} and extended to linear transport in \cite{mieussens2013}. A second order extension is suggested and validated. The generic nature of this method is also demonstrated in an application to the M2 model. Several test cases show the performances of this new scheme in both the M1 and M2 case.
\end{abstract}
\end{frontmatter}

\section{Introduction}

Kinetic equations appear in many fields of study such as plasma physics, radiative transfer, neutron transport and rarefied gas dynamics to model the dynamics of systems of particles. As the particle distribution is described in the phase space over time, accurately solving these equations  is expensive in terms of computational power. Furthermore macroscopic models correctly describe the system of particles as long as the Knudsen number (denoted by $\epsilon$)  remains low, which is defined as the ratio between the mean free path of the particles and a macroscopic length. The associated equations are much less costly to solve but their physical validity domain is limited. To describe transitional regimes and take into account kinetic effects without solving the complete equation, moments models are developed.

 These models aim to reduce the number of kinetic variables by closing a moment hierarchy of the kinetic equation. Closing the system consists in giving an expression of the unknown highest order moment as a function of the lower order ones. Such a relation can be provided by assuming the shape of the particle distribution at the microscopic scale in terms of the macroscopic variables. For example, the PN model rests on a Legendre series expansion of the distribution function under the small anisotropy hypothesis. The corresponding closure is linear, however the polynomial ansatz does not ensure the positivity of the distribution function \cite{modelcoll}.  In contrast, the MN model is based on the minimum entropy principle and guarantees this property for the Boltzmann entropy. Moreover, the MN system is hyperbolic and the flux limitation and entropy dissipation properties are ensured \cite{levermore1996moment,dubroca1999etude,alldredge2012high}. 

 Besides the high dimensional context, without specific treatment, numerical schemes for the kinetic equation or the moment model can be very expensive as they must resolve the smallest microscopic scale in the domain which constrains the space discretization and the time step for stability reasons. Furthermore, the limit scheme may not be consistent with the macroscopic model as the Knudsen number tends to zero. Asymptotic-preserving (AP) schemes have been developed to cope with this problem. Those schemes are consistent with the limit model and uniformly stable with $\epsilon$. They were first studied for neutron transport in \cite{larsen1987,larsen1989} and later in \cite{jin1991,jin1993}. In \cite{jin2000,klar1998}, AP schemes are obtained by decomposing the distribution function around the equilibrium and similar ideas are employed in \cite{buet2002diffusion,klar2001numerical,lemou2008new,bennoune2008uniformly,carrillo2008simulation,carrillo2008numerical}. Other approaches have been proposed in \cite{gosse2011transient} (well-balanced method) or \cite{lafitte2012asymptotic} (asymptotic-preserving projective integration scheme).

The Unified Gas Kinetic Scheme (UGKS) is an innovative AP scheme originally developed by Xu and Huang in 2010 in the context of rarefied gas dynamics \cite{kxu2010ugks}. Since then, it has been further improved and the general ideas have been applied to complex gas flows \cite{liu2017unified} (see \cite{zhu2021first} for other references). The UGKS was also extended to linear models with the diffusion limit in \cite{mieussens2013,sun2015radiativetransfer}.

For moment models, asymptotic-preserving schemes are usually constructed independently of the underlying kinetic equation. In most cases, a modification of the approximate Riemann solver is introduced to obtain the correct asymptotic behavior \cite{turpault,guisset2018,guisset2018admissible}. 

The main objective of this paper is to demonstrate how the UGKS may be utilized to develop a numerical scheme for the M1 moment model associated with a simple linear transport kinetic equation. Our idea is to apply the M1 closure at the numerical level on the numerical approximation (UGKS) of the linear kinetic equation. We prove that this scheme accurately captures the diffusion regime. Moreover, we suggest a second order extension that does not compromise the asymptotic-preserving property. Additionally, we show that the method developed in the M1 case is generic and can be applied to other moment models. In particular, a numerical scheme for the M2 moment model is given.

The outline of our article is as follows. First, in section \ref{sec_model} we briefly present the linear kinetic equation and the corresponding M1 model as well as their fundamental properties. Then, in section \ref{sec_ugks} the UGKS construction is summarized, the scheme for the M1 model is presented, and the second order extension is proposed. Next in section \ref{sec_m2}, a numerical scheme for the M2 moment model is also given. Finally, the schemes are validated in section \ref{sec_res}.

\section{The M1 closure for the linear transport }\label{sec_model}
\subsection{The linear transport equation}
The linear transport equation is a kinetic equation that describes the evolution of the particle number density $\phi$ as a function of time $t$, of space position $\mathbf r $ in $ \mathscr{D}$ an open set of $\mathbb{R}^3$ and of velocity direction $\mathbf \Omega $ in $ \mathscr{S}^2$ the unit sphere in 3 dimension space:
\begin{equation}
    \label{transfert_radiatif}
    \dfrac{1}{c}\partial_t \phi + \mathbf{\Omega} \cdot \nabla_{\mathbf{r}} \phi = \sigma (\dfrac{1}{4\pi}\int_{\mathscr{S}^2} \phi \mathrm{d} {{\Omega}} - \phi ).
\end{equation}
The number density represents the amount of particles in a given phase space volume at a certain time. From a physical point of view, this equation expresses the time variation of the number density through a collision operator in the absence of external forces. On the left-hand side, the total derivative in time describes the particles advection at velocity $c$ in the direction $\mathbf \Omega$. On the right-hand side, the collision operator models the particles interactions with the medium depending on the opacity $\sigma(\mathbf r)$ and reflects the rate of change of the number density. In this case, a linear relaxation operator is considered instead of the full non-linear Boltzmann one. This operator acts as a relaxation term towards the equilibrium state, which is the uniform velocity distribution. It preserves some basic fundamental properties such as mass conservation and entropy dissipation. 

In a small opacity medium, the particles are advected on the microscopic scale without colliding; this is the free transport regime. In that case, the number density is constant along the trajectories. Besides  the collision mechanism predominates and a global macroscopic diffusion behavior emerges when the opacity is high. 

 To study the diffusion regime and for computational purposes, it is convenient to work with the non-dimensional equation. In order to obtain this equation, several non-dimensional variables are introduced: $\phi'=\phi / \phi^*$, $t'=t / \tau$, $r'=r / L$, $\sigma'=\sigma / \sigma^*$ where $\tau$ is a characteristic time, $L$ a characteristic length and $\sigma^*$ a characteristic opacity homogeneous to the inverse of a length $\lambda$. This physical parameter represents the mean free path of a particle, that is, the average distance covered by a particle without a collision. Two non-dimensional numbers are introduced: the Knudsen number $\epsilon$ which is the ratio between the mean free path and the macroscopic length and $\eta$ which is the ratio between the macroscopic velocity and $c$: 
\begin{equation}
    \begin{matrix}
    \epsilon=\dfrac{\lambda}{L}, & \eta=\dfrac{L/\tau}{c}.
    \end{matrix}
\end{equation}
By omitting the prime symbol, the kinetic equation can be rewritten as a function of these quantities:
\begin{equation}
    \eta \partial_t \phi + \mathbf{\Omega} \cdot \nabla_{\mathbf{r}} \phi = \dfrac{\sigma}{\epsilon} (\dfrac{1}{4\pi} \sint{\mathscr{S}^2}{}{\phi}{\Omega}-\phi).
\end{equation}

In this article, we assume that the number density only depends on the slab axis variable $x$. In that case, the average of $\phi(t,\mathbf r, \Omega_x, \cdot, \cdot)$, denoted by $f(t,x,v)$ (where $v=\Omega_x$), satisfies the following one-dimensional equation:
\begin{equation}
    \label{equation_base}
 \partial_t f + \dfrac{v}{\eta} \partial_x f = \nu(\rho-f),
\end{equation}
where $\nu(x)=\frac{\sigma(x)}{\epsilon\eta}$ is the collision frequency and $\rho$ is the distribution function density: $\rho=\crochet{f}=\frac{1}{2}\int_{-1}^1f(\cdot,\cdot,v)\mathrm{d}v$. Integrating this kinetic equation over the velocity variable allows to get the macroscopic mass conservation equation:
\begin{equation}
    \label{conservation_masse}
    \epsilon \partial_t \rho + \partial_x j = 0,
\end{equation}
where $j=\crochet{vf}$ is the flux density.

\subsection{Asymptotic regimes}

As the Knudsen number $\epsilon$ tends to zero, the collision mechanism predominates at the microscopic scale and, as a consequence, the distribution function tends to its own density (at the first order in $\epsilon$). On the macroscopic scale, a global diffusion behavior emerges. To observe this phenomenon, the observation scale needs to coincide with the collision one, which implies $\eta=\epsilon$. In that case, a Hilbert expansion of the distribution function can be used to demonstrate that the density satisfies a diffusion equation at the first order in $\epsilon$:
\begin{equation}
    \label{equation_diffusion}
    \partial_t \rho = \partial_x \left(\kappa \partial_x \rho\right) + \mathcal{O}(\epsilon),
\end{equation}
where the diffusion coefficient is $\kappa(x)=\frac{1}{3\sigma(x)}$. Conversely, in the free transport regime, $\epsilon$ tends to infinity while $\eta$ remains constant. In that case, the limit equation is the usual linear advection equation without a source term:
\begin{equation}
    \label{equation_transport}
    \eta \partial_t f + v \partial_x f = 0.
\end{equation}
The particles are advected at their own speed $v/\eta$ without interacting with the medium.

\subsection{Entropy}
Due to the collision process, the particles tend to locally reach the equilibrium distribution (which is the uniform distribution) in a characteristic time $\tau=\frac{\eta \epsilon}{\sigma}$. From a physical point of view, a small perturbation out of that state leads to an increase of the physical entropy in the domain before returning to equilibrium. Mathematically, this irreversible process can be characterized by the local entropy inequality:
\begin{equation}
\label{h-theo-1}
    \eta \partial_t \crochet{g(f)} + \partial_x \crochet{v g(f)} \leq 0,
\end{equation}
where $g$ is any convex function. In a closed system, with suitable boundary conditions, the corresponding mathematical entropy $\mathscr{H}(t)=\int_\mathscr{D}\crochet{g(f(t,x,\cdot))}\mathrm{d}x$ is non-increasing;
\begin{equation}
\label{h-theo-2}
\derivee{\mathscr{H}}{t}(t)\leq 0.
\end{equation}
Numerical schemes for the kinetic equation should preserve this property which is a good indication of the system evolution.

\subsection{The M1 moment closure}
In a general context, solving kinetic equations is expensive due to the high dimensionality of the problem. In several physical applications, assumptions can be made on the shape of the distribution function. Thus, reduced models in velocity can be developed to lower the problem dimension and as a consequence the computational cost. A general procedure for elaborating such a model is to establish a moment hierarchy of the kinetic equation and then to choose a specific ansatz for the distribution function to close the resulting system. 

 The simplest hierarchy which enables the restoration of an angular anisotropy is obtained by integrating equation (\ref{equation_base}) against the vector $\mathbf{m}(v)=\begin{pmatrix} 1 & v \end{pmatrix}^T$ with respect to the velocity variable:
\begin{equation}
    \label{modele_moment}
    \partial_t\mathbf{U} + \partial_x \mathbf{F}(\mathbf{U}) = \nu \mathbf{S}(\mathbf{U}),
\end{equation}
where $\mathbf{U}=\begin{pmatrix} \rho & j\end{pmatrix}^T$ is the vector of conservative variables, $\mathbf{F}(\mathbf{U})=\frac{1}{\eta}\begin{pmatrix} j & q \end{pmatrix}^T$ is the flux vector where $q=\crochet{ v^2f}$ and $\mathbf{S}(\mathbf{U})=\begin{pmatrix} 0 & -j \end{pmatrix}^T$ is the source term. The first equation is the mass conservation equation (\ref{conservation_masse}). For any hierarchy, the integration process introduces a last unknown flux (in this case $q$) which can not be expressed, a priori, as a function of the previous moments. The M1 closure relies on an entropic argument to enforce the distribution function shape and to compute this flux as a function of the density and of the velocity $u=j/\rho$. Linked to this closure is the notion of moments realizability:
\begin{defi} [Moment realizability]
A moment vector $\mathbf{U}$ is realizable if there exists a non-negative distribution function $f$ such that $\crochet{\mathbf{m}f}=\mathbf{U}$.
\end{defi}
\begin{prop}
Let $\mathbf{U}=\begin{pmatrix} \rho & j\end{pmatrix}^T$ and $u=j/\rho$. The moment vector is realizable if and only if $\rho > 0$ and $|u|<1$, or $ \mathbf{U}=\mathbf{0}$.
\end{prop}
\begin{proof}
If $\mathbf{U}$ is realizable then $\rho=\crochet{f}\geq 0$. If $\rho=0$, then $f=0$ and hence $j=0$. If $\rho> 0$, then $|j|\leq \crochet{|v|f}< \crochet{f}=\rho$ since $|v|\leq 1$, and hence $|u|< 1$. The converse statement can be proven by setting $f=\hat{f}$ as defined in proposition \ref{prop:distribm1}. 
\end{proof}

\begin{prop}[M1 distribution function] \label{prop:distribm1}
Let $\mathbf{U}$ be a vector of realizable moments. If the density is non-zero, then the distribution function $\hat{f}$ which minimizes the Boltzmann entropy functional $h(f)= \crochet{f\ln{f}-f}$ under the constraint $\crochet{\mathbf{m} \hat{f}}=\mathbf{U}$ is
\begin{equation}
    \hat{f}(v)=e^{\mathbf\Lambda\cdot\mathbf{m}(v)}=\rho \dfrac{\beta}{\sinh{\beta}} e^{\beta v},
\end{equation}
where $\mathbf\Lambda = \begin{pmatrix} \alpha & \beta \end{pmatrix}^T$ is the vector of entropic variables and $\alpha=\ln{(\rho\frac{\beta}{\sinh\beta})}$.
The anisotropic variable $\beta$ is implicitly defined through the relation $u=z(\beta)$ where $z(\beta)= \coth{\beta}-\beta^{-1}$ is an invertible odd function in $[-1,1]$, continuously extendable at $\beta=0$.
\end{prop}
\begin{proof}
The M1 distribution function $\hat{f}$ satisfies the following constrained minimisation problem:
\begin{equation}
    \hat{f}=\argmin_{f\in\mathscr{S}}\crochet{f\ln{f}-f},
\end{equation}
where $\mathscr{S}=\{f\in \mathrm{L}^2([-1,1],\mathbb{R}_+)\text{ such that }\crochet{\mathbf{m}f}=\mathbf{U}\}$. The method of Lagrangian multipliers allows to show that:
\begin{equation}
    \hat{f}(v)=e^{\mathbf\Lambda \cdot \mathbf{m}(v)},
\end{equation}
where $\mathbf\Lambda \in \mathbb{R}^2$ is the Lagrangian multiplier vector. It can be implicitly expressed as a function of the conservative variable vector $\mathbf{U}$:
\begin{equation} \label{relation_multi}
    \begin{aligned}
    \mathbf{U}&=\crochet{\mathbf{m}e^{\mathbf\Lambda \cdot \mathbf{m}}} = \begin{pmatrix}
    \dfrac{e^\alpha}{\beta}\sinh{\beta} \\
    \dfrac{e^\alpha}{\beta}\sinh{\beta} \left(\coth\beta -\dfrac{1}{\beta}\right)
    \end{pmatrix}.\\
    \end{aligned}
\end{equation}
Thus, the M1 distribution function can be rewritten in terms of $\rho$ and $\beta$ and the relation between the anisotropic variable and the velocity appears.
\end{proof}
 Imposing the shape of the distribution function allows to close the system:
\begin{prop}[M1 closure]
The third moment $q$ of the M1 distribution function $\hat{f}$ is:
\begin{equation} \label{q_m1}
    q=\crochet{v^2\hat{f}}=\rho\left(1-2\dfrac{u}{\beta}\right).
\end{equation}
\end{prop}
System (\ref{modele_moment}) closed with relation (\ref{q_m1}) is the M1 model of the linear transport. We can notice that as the velocity tends to zero, $q$ tends to $\rho / 3$ which is nothing but the usual P1 closure. In the particular case of a zero density, the closing procedure is not applicable because the velocity and hence $\beta$ are not well defined anymore. But the continuity of $h$ at $f=0$ allows to set $\hat{f}=0$ and therefore $q=0$. \\

The following results hold on this model (see \cite{dubroca1999etude}).
\begin{prop}[System structure]
System (\ref{modele_moment})-(\ref{q_m1}) is hyperbolic (the Jacobian matrix of the system is diagonalizable and its eigenvalues are real) and ensures the moments realizability.
\end{prop}

\begin{prop}[Diffusion limit]
The density $\rho$ satisfies the diffusion equation $(\ref{equation_diffusion})$ at first order in $\epsilon$.
\end{prop}
 The validity domain of this model is directly linked to the quality of the distribution function projection on the set of M1 functions. As long as the distribution functions are close to this set, the model remains accurate. Two different types of distributions are well represented: the ones close to the equilibrium and the ones where the velocity is high. As soon as the distribution functions are far from the set of representable functions, this model becomes irrelevant.

\section{A UGKS based numerical scheme for the M1 model }\label{sec_ugks}
Developing a numerical scheme for the M1 hyperbolic system presents challenges for asymptotic preserving considerations. At first sight, a standard Riemann solver may appear suitable. However, without special treatment of the source term, it would not correctly capture the correct diffusion limit as the Knudsen number tends to zero.

 Several solvers rely on specific numerical fluxes designed to correctly capture the diffusion limit. For example in \cite{turpault,guisset2018}, the HLL approximate Riemann solver is modified by introducing a third stationary wave and by adjusting the nonlinear wave speed. As multiple choices are eligible to recover the correct asymptotic behavior, a particular attention is paid to the convergence speed to the diffusion regime as the Knudsen number tends to 0. 

An alternate and general procedure is to rely on a robust scheme for the kinetic equation. In this section an adaptation of the Unified Gas Kinetic Scheme (UGKS) for this model is explained.

\subsection{UGKS}
Since our new scheme is based on the UGKS, the solver construction for linear models with diffusion limit is adapted from \cite{mieussens2013} and summarized below.

\subsubsection{A finite volume formulation}
Let $[x_{i-1/2},x_{i+1/2}]$ be a control volume of size $\Delta x$ and $[t_n,t_{n+1}]$ be a time interval of size $\Delta t$. We define the averages of the density and distribution function on cell $i$ at time $t_n$ 
\begin{equation*} \begin{pmatrix}\rho_i^n \\ f_i^n(v) \end{pmatrix}=\dfrac{1}{\Delta x}\int_{x_{i-1/2}}^{x_{i+1/2}}\begin{pmatrix}\rho(t_n,x) \\ f(t_n,x,v)\end{pmatrix}\mathrm{d}x, \end{equation*}
and the macroscopic and microscopic numerical fluxes across the interface $x_{i+1/2}$
\begin{equation*} \begin{pmatrix} \Phi_{i+1/2} \\ \phi_{i+1/2}(v) \end{pmatrix}=\dfrac{1}{\eta \Delta t}\int_{t_n}^{t_{n+1}}\begin{pmatrix} \crochet{vf(t,x_{i+1/2},v)} \\ vf(t,x_{i+1/2},v) \end{pmatrix}\mathrm{d}t.\end{equation*}
The finite volume formulations of both the kinetic equation and the macroscopic conservation law are obtained by integrating equations ($\ref{equation_base}$)-($\ref{conservation_masse}$) over the control volume and over the time interval. These formulations emphasize the evolution of the volume averages through the cell interface fluxes between the two instants:
\begin{subequations}
\begin{align}
\dfrac{\rho_i^{n+1}-\rho_i^n}{\Delta t} + \dfrac{1}{\Delta x} (\Phi_{i+1/2}-\Phi_{i-1/2})&=0, \\
    \label{ugks_FV_micro} \dfrac{f_i^{n+1}-f_i^n}{\Delta t} + \dfrac{1}{\Delta x} (\phi_{i+1/2}-\phi_{i-1/2})&= \nu_i (\rho_i^{n+1} - f_i^{n+1}).
\end{align}
\end{subequations}
An implicit approximation of the collision term is chosen to obtain an asymptotically stable scheme. Developing a finite volume scheme for equation ($\ref{equation_base}$) involves giving a consistent and conservative approximation of the microscopic numerical flux $\phi_{i+1/2}$ and therefore of the macroscopic one $\Phi_{i+1/2}=\crochet{\phi_{i+1/2}}$. At this stage, the velocity variable $v$ is kept continuous and omitted.

\subsubsection{A characteristic based approach}
The main idea of UGKS is to rely on the integral representation of the kinetic equation solution (given by the method of characteristics) to elaborate the numerical flux. This way, the collision term is naturally taken into account. In case of constant opacity (\ref{equation_base}) is equivalent to:
\begin{equation}
\label{representation_inte}
    \derivee{}{t}\left(e^{\nu t} f(t,x+\dfrac{v}{\eta}t,v)\right)=\nu e^{\nu t}\rho(t,x+\dfrac{v}{\eta}t).
\end{equation}
Assuming the opacity variations are negligible at the scale of a cell and of a time step, we consider this expression as an approximation around each cell. Relation (\ref{representation_inte}) is then evaluated at the interface $x_{i+1/2}$ and integrated between two given times, $t_n$ and $t>t_n$, which gives
\begin{equation}
    \label{representation_interface}
    \begin{aligned}
    f(t,x_{i+1/2},v)&\approx e^{-\nu_{i+1/2}(t-t_n)} f(t_n,x_{i+1/2}-\dfrac{v}{\eta}(t-t_n),v) \\
    ~&+\nu_{i+1/2}\int_{t_n}^t  e^{-\nu_{i+1/2}(t-s)} \rho(s,x_{i+1/2} - \dfrac{v}{\eta}(t-s)) \mathrm{d}s,
    \end{aligned}
\end{equation}
where $\nu_{i+1/2}={\sigma_{i+1/2}}/{\eta\epsilon}$ is the collision frequency at the interface. The total number of particles at the interface can be separated into two categories ; advected and scattered particles. Depending of the collision frequency, some particles do not interact with others and are simply transported from the foot of the characteristic $x_{i+1/2}-\frac{v}{\eta}(t-t_n)$ to the interface. Other particles have a certain probability of colliding once at some time $s$ such as $t>s>t_n$ and acquiring the specific $v$ velocity at $x_{i+1/2}-\frac{v}{\eta}(t-s)$. All of these particles are then transported to the interface. 

 In order to evaluate the numerical flux from relation ($\ref{representation_interface}$), distribution function and density reconstructions in space and time need to be introduced. Appropriate choices are mandatory to preserve the asymptotics and achieve second order convergence in space. The reconstructions are:
\begin{subequations}
\begin{equation}
    \label{reconstruction_densite}
    \rho(t,x)=\begin{cases}
    \rho_{i+1/2}^n + \delta_x^L \rho_{i+1/2}^{n} (x-x_{i+1/2}) & \text{if } x < x_{i+1/2} \\
    \rho_{i+1/2}^n + \delta_x^R \rho_{i+1/2}^{n} (x-x_{i+1/2}) & \text{if } x > x_{i+1/2} \\
    \end{cases},
\end{equation}
\begin{equation}
    \label{reconstruction_f}
    f(t_n,x,v)= \begin{cases}
    f_i^n  + \delta_x f_i^n (x-x_i) & \text{if } x < x_{i+1/2} \\
    f_{i+1}^{n} + \delta_x f_{i+1}^n (x-x_{i+1}) & \text{if } x > x_{i+1/2}
    \end{cases},
\end{equation}
\end{subequations}
where $\delta_x^{LR} \rho_{i+1/2}^{n}$ are the left and right finite differences slopes:
\begin{equation}
    \begin{matrix}
    \delta_x^L \rho_{i+1/2}^{n} = \dfrac{\rho_{i+1/2}^n-\rho_i^n}{\Delta x / 2}, &
    \delta_x^R \rho_{i+1/2}^{n} = \dfrac{\rho_{i+1}^n-\rho_{i+1/2}^n}{\Delta x / 2},
    \end{matrix}
\end{equation}
and the interface density $\rho_{i+1/2}$ is:
\begin{equation}
    \rho_{i+1/2}=\crochet{\indi{f_i^n}{f_{i+1}^n}{v}}=\rho_i^{n+}+\rho_{i+1}^{n-}.
\end{equation}
The choice for that density turns out to be not that important. For example the mean value $\frac{\rho_i+\rho_{i+1}}{2}$ may be appropriate. However, from a physical point a view, using the half densities on each side seems to be equally relevant in order to consider the real distribution of the density near the interface and to ensure the BGK compatibility condition at $t=t_n$ in (\ref{representation_interface}):
\begin{equation*}
    \crochet{\rho-f}(t,x_{i+1/2})=0.
\end{equation*}
The distribution function slopes need to be limited to ensure the decrease of the total variation. Let $\psi$ be a TVD slope limiter (for example the van Leer limiter is given by $\psi(x,y)=(\sign(x)+\sign(y)) \frac{|x||y|}{|x|+|y|}$ ). Then the slope is given by:
\begin{equation}
    \delta_x f_i^n = \psi \left( \dfrac{f_{i+1}^n-f_{i}^n}{\Delta x},\dfrac{f_{i}^n-f_{i-1}^n}{\Delta x} \right).
\end{equation}
 To evaluate the numerical flux $\phi_{i+1/2}$, the reconstructed quantities are employed in ($\ref{representation_interface}$) before time integration. It should be noted that in the diffusion limit, the foot of the characteristics might be arbitrarily far from the interface. However, due to the collision mechanism, the particles are constrained near the interface (as shown by the exponential term in ($\ref{representation_interface}$)). Therefore, it is legitimate to neglect the influence of remote particles by extending the reconstructions validity domain. Finally, the microscopic numerical flux takes the following form
\begin{equation}
\label{flux_micro_2}
    \begin{aligned}
    \phi_{i+1/2}(v) = & A_{i+1/2} v \left(f_i^{n(+)} \mathbb{1} _{v>0} + f_{i+1}^{n(-)} \mathbb{1}_{v<0}\right) \\
     + &B_{i+1/2} v^2 (\delta_x f_{i}^n \mathbb{1}_{v>0}+\delta_x f_{i+1}^n \mathbb{1}_{v<0})\\
     + &C_{i+1/2} v \rho_{i+1/2}^n \\
     + &D_{i+1/2} v^2 (\delta_x^L \rho_{i+1/2}^{n} \mathbb{1} _{v>0} + \delta_x^R \rho_{i+1/2}^{n} \mathbb{1} _{v<0}),
    \end{aligned}
\end{equation}
and the macroscopic one is
\begin{equation}
\label{flux_macro_2}
\begin{aligned}
\Phi_{i+ 1/2}= & A_{i+ 1/2} \crochet{v f_i^{n(+)} \mathbb{1} _{v>0} +v f_{i+1}^{n(-)} \mathbb{1}_{v<0}}  \\
+ & B_{i+1/2} \crochet{v^2\delta_x f_i^n \mathbb{1} _{v>0} + v^2\delta_x f_{i+1}^n \mathbb{1}_{v<0}}\\
+ & \dfrac{D_{i+1/2}}{3\Delta x} (\rho_{i+1}^{n}-\rho_i^n),
\end{aligned}
\end{equation}
where $f_i^{n(\pm)}=f_i^n\pm\frac{\Delta x}{2}\delta_xf_i^n$. The integration coefficients $A_{i+ 1/2}$, $B_{i+ 1/2}$, $C_{i+ 1/2}$ and $D_{i+ 1/2}$ are interface values of functions
\begin{subequations}
\begin{align}
A(\Delta t,\eta,\epsilon,\sigma)&=\dfrac{-1}{\eta}\dfrac{(1-e^w)}{w}, \\
B(\Delta t,\eta,\epsilon,\sigma)&=\dfrac{1}{\sigma} \dfrac{\epsilon}{\eta} \left( e^{w} +\dfrac{1-e^w}{w} \right), \\
C(\Delta t,\eta,\epsilon,\sigma)&=\dfrac{1}{\eta}\left(1+\dfrac{1-e^w}{w}\right),\\    
D(\Delta t,\eta,\epsilon,\sigma)&=\dfrac{-1}{\sigma} \dfrac{\epsilon}{\eta}\left(1+e^{w}+2\dfrac{1-e^w}{w}\right),
\end{align}
\end{subequations}
at $\sigma_{i+1/2}=\frac{\sigma_i+\sigma_{i+1}}{2}$ and where $w=-\nu \Delta t$.

\subsubsection{Asymptotic behaviour and stability}
We examine the asymptotic preserving property of the scheme in both the diffusion and free transport regimes. The opacity $\sigma$ is assumed bounded in the asymptotic analysis. The numerical fluxes behavior is entirely determined by the integration coefficients limits. In the diffusion limit, the constraint $\eta=\epsilon$ is enforced, the following limits hold:
\begin{equation*}
    \begin{matrix}
    A(\Delta t,\epsilon,\epsilon,\sigma)\underset{\epsilon \to 0}{\longrightarrow}0, &
    B(\Delta t,\epsilon,\epsilon,\sigma)\underset{\epsilon \to 0}{\longrightarrow}0, &
    C(\Delta t,\epsilon,\epsilon,\sigma)\underset{\epsilon \to 0}{\sim}\dfrac{1}{\epsilon}, &
    D(\Delta t,\epsilon,\epsilon,\sigma)\underset{\epsilon \to 0}{\longrightarrow}\dfrac{-1}{\sigma}. &
    \end{matrix}
\end{equation*}
As a consequence, the limit macroscopic flux is:
\begin{equation}
    \label{flux_macro_limit}
    \Phi_{i+1/2}\underset{\epsilon \to 0}{\longrightarrow} \frac{-1}{3 \sigma_{i+1/2}} \dfrac{\rho_{i+1}^n-\rho_i^n}{\Delta x},
\end{equation}
which is the usual second order flux for the diffusion equation. The correct diffusion coefficient $\kappa(x)=\frac{1}{3\sigma(x)}$ is recovered. In the free transport regime, obtained with constant $\eta$ and large $\epsilon$, we have the following limits:
\begin{equation*}
    \begin{matrix}
    A(\Delta t,\eta,\epsilon,\sigma)\underset{\epsilon \to \infty}{\longrightarrow}\dfrac{1}{\eta}, &
    B(\Delta t,\eta,\epsilon,\sigma)\underset{\epsilon \to \infty}{\longrightarrow}\dfrac{-\Delta t}{2\eta^2}, &
    C(\Delta t,\eta,\epsilon,\sigma)\underset{\epsilon \to \infty}{\longrightarrow}0,&
    D(\Delta t,\eta,\epsilon,\sigma)\underset{\epsilon \to \infty}{\longrightarrow}0.&
    \end{matrix}
\end{equation*}
The limit microscopic flux is:
\begin{equation}
    \label{flux_micro_limit}
    \phi_{i+1/2} \underset{\epsilon \to \infty}{\longrightarrow} \dfrac{v}{\eta} (f_i^{n(+)} \mathbb{1}_{v>0} + f_{i+1}^{n(-)} \mathbb{1}_{v<0}) 
    -\Delta t \dfrac{v^2}{2\eta^2} (\delta_x f_i^n \mathbb{1}_{v>0} + \delta_x f_{i+1}^n \mathbb{1}_{v<0}),
\end{equation}
which is a second order in space and time flux for the free transport equation. Even if we are not able to mathematically prove that this scheme is uniformly stable in some sense under a CFL condition, it is observed that the following heuristic condition is sufficient (see \cite{mieussens2013}):
\begin{equation}
    \label{cfl_somme}
    \Delta t \leq \dfrac{3}{2}\sigma\Delta x^2  + \eta \Delta x.
\end{equation}

\subsection{UGKS-M1}
\subsubsection{An entropic closure of the UGKS}
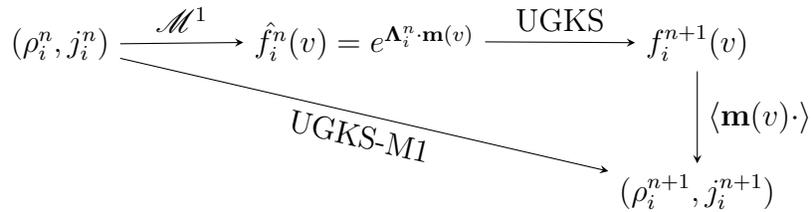
\begin{figure}[!htbp]
    \centering
    \begin{tikzpicture}
      \matrix (m) [matrix of math nodes,row sep=3em,column sep=4em,minimum width=2em]
      {
         (\rho_i^n,j_i^n) & \hat{f_i^n}(v)=e^{\boldsymbol\Lambda_i^n \cdot \textbf{m}(v)} & f_i^{n+1}(v) \\
          ~ & ~ &  (\rho_i^{n+1},j_i^{n+1}) \\
          };
      \path[-stealth]
        (m-1-1) edge node [above] {$\mathscr{M}^1$} (m-1-2)
        (m-1-2) edge node [above] {UGKS} (m-1-3)
        (m-1-1) edge node [midway, sloped, below] {UGKS-M1} (m-2-3)
        (m-1-3) edge node [right] {$\crochet{\mathbf{m}(v) \cdot }$} (m-2-3);
    \end{tikzpicture}
    \caption{Structure of the UGKS-M1 scheme \label{fig:ugks_m1}}
\end{figure}

The natural idea behind this new scheme is to apply the UGKS to the M1 distribution function $(\hat{f}_i^n)$ reconstructed from the moments $(\rho_i^n,j_i^n)$. Then, the macroscopic variables at time $t_{n+1}$ are the moments of  $(f_i^{n+1})$. This process is globally represented in figure \ref{fig:ugks_m1}. From another point of view, this procedure can be seen as a systematic projection of the distribution function in the M1 set at each time step in UGKS. This new scheme then appears as a M1 moment closure of UGKS. 

Taking the first two moments of the microscopic scheme (\ref{ugks_FV_micro}) provides a finite volume formulation for the vector of discrete conservative variable $\mathbf{U}_i^n=\begin{pmatrix} \rho_i^n & j_i^n \end{pmatrix}^T$:

\begin{equation}\label{ugsk_m1_fv}
    \dfrac{\mathbf{U}_i^{n+1}-\mathbf{U}_i^{n}}{\Delta t} + \dfrac{1}{\Delta x}(\mathbf{\Phi}_{i+1/2}-\mathbf{\Phi}_{i-1/2})=\nu_i\mathbf{S}(\mathbf{U}_i^{n+1}),
\end{equation}
where $ \mathbf{\Phi}_{i+1/2}=\crochet{\mathbf{m}(v)\phi_{i+1/2}(v)}=\begin{pmatrix} \Phi_{i+1/2}^\rho & \Phi_{i+1/2}^j \end{pmatrix}^T$. Then, this macroscopic flux vector is computed by integrating the microscopic UGKS flux (\ref{flux_micro_2}) with $f_i^n=\hat{f}_i^n$. First, we define the fluxes without the second order term in (\ref{reconstruction_f}) (the distribution function reconstruction is constant per cell):
\begin{subequations}
    \label{ugks_m1_flux}
    \begin{align}\label{ugks_m1_flux_p}
    {\Phi}_{i+1/2}^\rho &= A_{i+ 1/2} \crochet{v\hat{f}_i^{n+} \mathbb{1} _{v>0}+ v\hat{f}_{i+1}^{n-} \mathbb{1}_{v<0}}
    + \dfrac{D_{i+1/2}}{3\Delta x} (\rho_{i+1}^{n}-\rho_i^n), \\
    {\Phi}_{i+1/2}^j &=A_{i+ 1/2} \crochet{v^2\hat{f}_i^{n} \mathbb{1} _{v>0}+ v^2\hat{f}_{i+1}^{n} \mathbb{1}_{v<0}}
    +\dfrac{C_{i+1/2}}{3} \rho_{i+1/2}^n.
    \end{align}
\end{subequations}
where $(\alpha_i,\beta_i)$ are the entropic variables associated with $\mathbf{U}_i^{n}$ as defined in proposition \ref{prop:distribm1}:
\begin{equation*}
\begin{matrix}
    e^{\alpha_i}=\rho_i^n \dfrac{\beta_i}{\sinh{\beta_i}}, &
    \beta_i=z^{-1}(\dfrac{j_i^n}{\rho_i^n}),
\end{matrix}
\end{equation*}
where $z^{-1}$ is the inverse function of $z(\beta)=\coth{\beta}-\beta^{-1}$. The anisotropic factor $\beta_i$ can be numerically computed using the Newton method. By nature, this new scheme is asymptotic preserving. Indeed, by performing the same analysis as with UGKS, we can notice that the first macroscopic flux (\ref{ugks_m1_flux_p}) tends to the correct diffusion flux (\ref{flux_macro_limit}) in the corresponding limit. \\

\noindent Three numerical difficulties appear: \begin{itemize}
    \item The exact value of the integrals of the form $\crochet{v^i \hat{f} \mathbb{1}_{v \lessgtr 0}}$ should be programmed in a developed form (see \ref{ann1}) with $e^\beta$ in factor to avoid an accumulation of round-off error at low velocities.
    \item Below a certain $\beta$ (or $u$) threshold, the same integrals $\crochet{v^i \hat{f} \mathbb{1}_{v \lessgtr 0}}$ should be set to the correct limit, which is $\rho_i^n \crochet{v^i \mathbb{1}_{v \lessgtr 0}}$.
    \item For low densities, $u_i^n=\frac{j_i^n}{\rho_i^n}$ may not be well-defined anymore. Below a certain density threshold, we set $\hat{f}_i^n=0$ and therefore $\beta_i=0$ to correctly compute the flux in this limit.
\end{itemize}

\subsubsection{Definition and realizability of the scheme}
A question addressed here is the definition of the scheme (\ref{ugsk_m1_fv})-(\ref{ugks_m1_flux}). It is clear that $\mathbf{U}_i^{n+1}$ can be computed only if $\hat{f}_i^n$ can be defined in every cell. This requires the moment vector $\mathbf{U}_i^{n}$ to be realizable. The scheme can only be iterated only if that property holds for every time step. In other words, we should prove that: $(\mathbf{U}_i^n$) is realizable implies that $(\mathbf{U}_i^{n+1})$ is realizable as well. Such a scheme is said to be realizable and this property could be obtained with the following simple argument (see \cite{desjardins2008quadrature}). Our scheme can be written in the form
\begin{equation}
    \mathbf{U}_i^{n+1}=\crochet{\mathbf{m}f_i^{n+1}},
\end{equation}
where $f_i^{n+1}$ is obtained with one time step of UGKS initialized with $\hat{f}_i^n$. Then, $\textbf{U}_i^{n+1}$ is realizable if $f_i^{n+1}$ is non-negative. Consequently, the realizability of the scheme can be reduced to the question of the positivity of UGKS, at least with initial data given by a M1 distribution. The proof of this property is not fully completed at the moment. However, from a practical point of view the reliability of the scheme seems to be ensured for smooth enough initial and boundary data and under the cfl-like condition (\ref{cfl_somme}).

\subsubsection{Second order in space}
In the previous part we dropped the linear part of the distribution function reconstruction in  (\ref{reconstruction_f}). This term is problematic as the integrals of the form $\crochet{v^2 \delta_x \hat{f} \mathbb{1}_{v \lessgtr 0}}$ cannot be analytically expressed as a function of the entropic variables due to the non-linearity introduced by the slope limiter. To achieve a second order convergence rate in space, a different reconstruction of the distribution function is used as proposed in \cite{kxu2010ugks} for the Boltzmann equation of rarefied gas dynamics. First, the vector of conservative variables is reconstructed:

\begin{equation}
    \mathbf{U}_i^n (x)=\begin{cases}
    \mathbf{U}_i^n + \delta \mathbf{U}_{i}^n(x-x_{i}) &\text{ if } x<x_{i+1/2}\\
    \mathbf{U}_{i+1}^n + \delta \mathbf{U}_{i+1}^n(x-x_{i+1}) &\text{ if } x>x_{i+1/2}
    \end{cases},
\end{equation}
where the finite difference slope is $\delta \mathbf{U}_i^n=\frac{1}{\Delta x}(\mathbf{U}_{i+1}^n-\mathbf{U}_i^n) \mathbf{\phi}(\mathbf{r}_i)$, $\phi$ is a slope limiter and $\mathbf{r}_i=\begin{pmatrix} \frac{\rho_i-\rho_{i-1}}{\rho_{i+1}-\rho_i} & \frac{j_i-j_{i-1}}{j_{i+1}-j_i} \end{pmatrix}^T$ is the local slope defined component-wise. Then, we expand $\hat{f}(\mathbf{U}_i^n(x))=\exp{(\mathbf{\Lambda}(\mathbf{U}_i^n(x))\cdot \mathbf{m})}$ in Taylor series (for example when $x<x_{i+1/2}$):
\begin{equation}
    \begin{aligned}
    \hat{f}(\mathbf{U}_i^n (x))&=\hat{f}(\mathbf{U}_i^n) + \dfrac{\mathbf{\mathrm{d}\hat{f}}}{\mathrm{d}\mathbf{U}}(\mathbf{U}_i^n) \cdot \mathbf{\delta U}_i^n (x-x_i) \\
    &=\hat{f}(\mathbf{U}_i^n) + \mathbf{J}_\mathbf\Lambda (\mathbf{U}_i^n)^T \mathbf{m}\hat{f}(\mathbf{U}_i^n)    \cdot \mathbf{\delta U}_i^n (x-x_i),
    \end{aligned}
\end{equation}
where the Jacobian matrix $\mathbf{J}_\mathbf\Lambda (\mathbf{U})$ is
\begin{equation}
\begin{aligned}
    \mathbf{J}_\mathbf\Lambda (\mathbf{U})=\mathbf{J}_\mathbf{U} (\mathbf{\Lambda})^{-1}&=\crochet{\mathbf m \otimes \mathbf m \exp{(\mathbf\Lambda \cdot \mathbf m )}}^{-1}, \\
    &=\dfrac{\rho^{-1}}{1-2\dfrac{u}{\beta}-u^2}
     \begin{pmatrix}
    1-2\dfrac{u}{\beta} & -u \\
    -u & 1 \\
    \end{pmatrix}.
    \end{aligned}
\end{equation}
Then, the M1 distribution function reconstruction is:
\begin{equation}
    \label{reconstruction_fm1}
    f(t_n,x,v)= \begin{cases}
    \hat{f}_i^n  + \delta_x \hat{f}(\mathbf{U}_i^n) (x-x_i) & \text{if } x < x_{i+1/2} \\
    \hat{f}_{i+1}^{n} + \delta_x \hat{f}(\mathbf{U}_{i+1}^n) (x-x_{i+1}) & \text{if } x > x_{i+1/2}
    \end{cases},
\end{equation}
where the slope is $\delta_x\hat{f}(\mathbf{U}_i^n)= \mathbf{J}_\mathbf\Lambda(\mathbf{U}_i^n) \mathbf{\delta U}_i^n \cdot \mathbf{m}\hat{f}(\mathbf{U}_i^n) $. Finally, the second order fluxes are:

\begin{subequations} \label{ugks_m1_flux2}
    \begin{align}
    \begin{split}
         {{\Phi}_{i+1/2}^\rho} &= A_{i+ 1/2} \crochet{v\hat{f}_i^{n+} \mathbb{1} _{v>0}+ v\hat{f}_{i+1}^{n-} \mathbb{1}_{v<0}} + B_{i+1/2}\crochet{v^2 \delta_x \hat{f}(\mathbf{U}_i^n)\mathbb{1}_{v>0} + v^2 \delta_x \hat{f}(\mathbf{U}_{i+1}^n)\mathbb{1}_{v<0}} \\
         &+ \dfrac{D_{i+1/2}}{3\Delta x} (\rho_{i+1}^{n}-\rho_i^n), \\
    \end{split} \\
    \begin{split}
        {{\Phi}_{i+1/2}^j} &=A_{i+ 1/2} \crochet{v^2\hat{f}_i^{n+} \mathbb{1} _{v>0} + v^2\hat{f}_{i+1}^{n-} \mathbb{1}_{v<0}} + B_{i+1/2}\crochet{v^3 \delta_x \hat{f}(\mathbf{U}_i^n) \mathbb{1} _{v>0} + v^3 \delta_x \hat{f}(\mathbf{U}_{i+1}^n) \mathbb{1} _{v<0}} \\
        &+ \dfrac{C_{i+1/2}}{3} \crochet{\hat{f}_i^n \mathbb{1} _{v>0} + \hat{f}_{i+1}^n \mathbb{1} _{v<0}}, \\
    \end{split}
    \end{align}
\end{subequations}
and where
\begin{equation*}
    \hat{f}_i^{n \pm}=\hat{f}_i^n \pm \dfrac{\Delta x}{2}\delta_x \hat{f}_i^n.
\end{equation*}
The full expressions of the fluxes can be found in \ref{ann1}. 

\section{Extension to other moment closures}\label{sec_m2}
\subsection{General Framework}
The method described in this article to obtain a numerical scheme for the M1 model is generic and can be easily applied to other moment models. Let $\mathbf{m}(v) \in \mathbb{R}^{d+1}$ be a vector composed of the elements of any basis of $\mathbb{R}_{d}[X]$ (the set of polynomials of degree at most $d$) and $\mathbf{U}=\crochet{\mathbf{m}f} \in \mathbb{R}^{d+1}$ be the corresponding moment vector of $f$. A moment model is obtained by approximating $f$ by some specific ansatz $\hat{f}(\mathbf{U})$ that realizes the same moments. The resulting moment model takes the same form as equation (\ref{modele_moment}):
\begin{equation}
    \label{modele_momentM2}
    \partial_t\mathbf{U} + \partial_x \mathbf{F}(\mathbf{U}) = \nu \mathbf{S}(\mathbf{U}),
\end{equation}
where $\mathbf{F}(\mathbf{U})=\frac{1}{\eta}\crochet{v\mathbf{m}\hat{f}(\mathbf{U})}$ is the flux vector and $\mathbf{S}(\mathbf{U})=\crochet{\mathbf{m}}\rho-\mathbf{U}$ is the source term. 

The finite volume scheme can be obtained as outlined in section 3, resulting in:
\begin{equation}\label{ugks_m2_fv}
    \dfrac{\mathbf{U}_i^{n+1}-\mathbf{U}_i^{n}}{\Delta t} + \dfrac{1}{\Delta x}(\mathbf{\Phi}_{i+1/2}-\mathbf{\Phi}_{i-1/2})=\nu_i\mathbf{S}(\mathbf{U}_i^{n+1}),
\end{equation}
where the numerical flux $\mathbf{\Phi}_{i+1/2}$ is
\begin{equation}
    \boldsymbol\Phi_{i+1/2}=\crochet{\mathbf{m}\phi_{i+1/2}},
\end{equation}
where $\phi_{i+1/2}$ is the UGKS microscopic flux (\ref{flux_micro_2}) in which the distribution function is set to the chosen ansatz: $f_i^n=\hat{f}(\mathbf{U}_i^n)$.

\subsection{Application to M2}
The M2 moment model is the next order (after M1) in the hierarchy of entropic moment models (see \cite{dubroca1999etude,hauck2011high,pichard2017approximation} for more details) given by $\textbf{m}(v)=\begin{pmatrix}1 & v & v^2\end{pmatrix}^T$. The moment vector is denoted by $\mathbf{U}=\begin{pmatrix} \rho & j & q\end{pmatrix}^T$. Similar to M1, the closure ansatz is obtained by minimizing the Boltzmann entropy to obtain
\begin{equation}
    \hat{f}(\mathbf{U})(v)=e^{\boldsymbol\Lambda(\mathbf{U}) \cdot \mathbf{m}(v)}=e^{\alpha + \beta  v + \gamma v^2},
\end{equation}
where $\boldsymbol\Lambda=\begin{pmatrix}\alpha & \beta & \gamma \end{pmatrix}^T$ is the vector of entropic variables. The M2 moment model then reads as in (\ref{modele_momentM2}) with the source term $\mathbf{S}(\mathbf{U})=\begin{pmatrix}0 & -j & \frac{\rho}{3}-q\end{pmatrix}^T$. This resulting M2 moment model can be shown to be a hyperbolic system (see \cite{dubroca1999etude}).

The general framework presented before leads to the UGKS-M2 scheme (\ref{ugks_m2_fv}), where after some algebra the numerical flux $\mathbf{\Phi}_{i+1/2}=\begin{pmatrix}\Phi_{i+1/2}^\rho & \Phi_{i+1/2}^j & \Phi_{i+1/2}^q \end{pmatrix}$ simplifies to
\begin{subequations}
    \begin{align}
            {{\Phi}_{i+1/2}^\rho} &= A_{i+1/2}\crochet{v\hat{f}_{i}^n\mathbb{1}_{v>0}+v\hat{f}_{i+1}^n\mathbb{1}_{v<0}}+ \dfrac{D_{i+1/2}}{3\Delta x} (\rho_{i+1}^{n}-\rho_i^n), \\
            {{\Phi}_{i+1/2}^j} &=A_{i+1/2}\crochet{v^2\hat{f}_{i}^n\mathbb{1}_{v>0}+v^2\hat{f}_{i+1}^n\mathbb{1}_{v<0}}+ \dfrac{C_{i+1/2}}{3}\rho_{i+1/2}^n,\\
            {{\Phi}_{i+1/2}^q} &=A_{i+1/2}\crochet{v^3\hat{f}_{i}^n\mathbb{1}_{v>0}+v^3\hat{f}_{i+1}^n\mathbb{1}_{v<0}}+ \dfrac{D_{i+1/2}}{5\Delta x} (\rho_{i+1}^{n}-\rho_i^n),
    \end{align}
\end{subequations}
with the half moments
\begin{subequations} \label{half_moments_m2}
\begin{align} \label{half_densities_m2}
    \rho^\pm &= \crochet{\hat{f} \mathbb{1}_{v\gtrless0}}= \dfrac{\pm1}{2\sqrt{|\gamma|}}e^\alpha 
    \begin{cases}
    e^{\pm\beta+\gamma}\Dp\left(\pm\sqrt{|\gamma|}+\dfrac{\beta}{2\sqrt{|\gamma|}}\right)-\Dp\left(\dfrac{\beta}{2\sqrt{|\gamma|}}\right) & \gamma > 0\\
    \Dm\left(\dfrac{-\beta}{2\sqrt{|\gamma|}}\right) -e^{\pm\beta+\gamma}\Dm\left(\pm\sqrt{|\gamma|}-\dfrac{\beta}{2\sqrt{|\gamma|}}\right) & \gamma < 0\\
    \end{cases}, \\
    j^\pm &= \crochet{v\hat{f} \mathbb{1}_{v\gtrless0}}=\dfrac{\pm1}{4\gamma}e^\alpha (e^{\pm\beta+\gamma}-1) - \dfrac{\beta}{2\gamma} \rho^\pm, \\
    q^\pm&= \crochet{v^2\hat{f} \mathbb{1}_{v\gtrless0}}=\dfrac{1}{4\gamma}e^{\alpha\pm\beta+\gamma} -\dfrac{\beta}{2\gamma}j^\pm -\dfrac{1}{2\gamma}\rho^\pm,\\
    k^\pm&= \crochet{v^3\hat{f} \mathbb{1}_{v\gtrless0}}=\dfrac{\pm1}{4\gamma}\dfrac{\gamma-1}{\gamma}e^\alpha (e^{\pm\beta+\gamma}-1)  -\dfrac{\beta}{2\gamma}q^\pm+\dfrac{\beta}{2\gamma^2}\rho^\pm,
\end{align}
\end{subequations}
and where $\Dp$ is the Dawson function and $\Dm$ is the scaled complementary error function (up to a constant factor) defined for every real number $x$ by
\begin{equation*}
\begin{aligned}
\Dp(x)=e^{- x^2}\int_0^x e^{ t^2}\mathrm{d}t \quad \text{ and } \quad
\Dm(x)=e^{x^2}\int_x^{+\infty} e^{ -t^2}\mathrm{d}t.
\end{aligned}
\end{equation*}
The parameters $\alpha,\beta,\gamma$ in (\ref{half_moments_m2}) are the components of the entropic variable $\boldsymbol\Lambda$, which can be computed by inverting the relation
\begin{equation} \label{realisation_moment}
\mathbf{U}=\crochet{\mathbf{m}\hat{f}(\boldsymbol\Lambda)}.
\end{equation}
This inversion must be performed numerically and is rather delicate. First, we found that for numerical reasons, it is relevant to write (\ref{realisation_moment}) as
\begin{subequations} 
    \begin{align}
    \rho&=\dfrac{1}{2\sqrt{|\gamma|}}e^{\alpha+\gamma} 
    \begin{cases}
    e^\beta \Dp\left(\sqrt{|\gamma|}+\dfrac{\beta}{2\sqrt{|\gamma|}}\right)-e^{-\beta} \Dp\left(-\sqrt{|\gamma|}+\dfrac{\beta}{2\sqrt{|\gamma|}}\right) & \gamma>0 \\
    e^{-\beta} \Dm\left(-\sqrt{|\gamma|}-\dfrac{\beta}{2\sqrt{|\gamma|}}\right)-e^{\beta} \Dm\left(\sqrt{|\gamma|}-\dfrac{\beta}{2\sqrt{|\gamma|}}\right) & \gamma<0 \\
    \end{cases} \label{p_m2}\\
    j&=\dfrac{1}{2\gamma}e^{\alpha+\gamma}\sinh\beta - \dfrac{\beta}{2\gamma} \rho,\label{j_m2} \\
    q&=\dfrac{1}{2\gamma}e^{\alpha+\gamma}\cosh\beta  -\dfrac{\beta}{2\gamma}j -\dfrac{1}{2\gamma}\rho.\label{q_m2}
    \end{align}
\end{subequations}
One common method to solve (\ref{realisation_moment}) for $\boldsymbol\Lambda$ is the Newton method. However, we found it to be insufficiently robust in many cases. Instead, note that (\ref{realisation_moment}) can be expressed as $\nabla J(\boldsymbol\Lambda)=0$ where the functional $J$ is $J(\boldsymbol\Lambda)=\crochet{\exp{(\boldsymbol\Lambda\cdot \mathbf{m})}} - \boldsymbol\Lambda \cdot \mathbf{U}$. This shows that $\boldsymbol\Lambda$ can be computed as the minimum of $J$ on $\mathbb{R}^{3}$, which is achieved with a standard gradient descent algorithm that turns out to be very robust. Note that similarly to the M1 model, the first entropic variable $\alpha$ can be eliminated and replaced by the density in (\ref{j_m2})-(\ref{q_m2}) using (\ref{p_m2}). 
\begin{rema}
The resulting UGKS-M2 scheme is not diagonally implicit due to the density term in the source term of the third equation. However, the density can be computed first before updating $q$. Alternatively, the scheme can be made diagonally implicit by using the Legendre basis $\mathbf{m}(v)=\begin{pmatrix}1 & v & \frac{3}{2}(v^2-\frac{1}{3})\end{pmatrix}^T$ to form the moment hierarchy.
\end{rema}
\begin{rema}
As $\gamma$ tends to $0$, the M2 distribution function tends to the M1 distribution function. However, the density expressions (\ref{p_m2}) and (\ref{half_densities_m2}) become singular in this limit. For numerical robustness, it is necessary to use asymptotic expansions of the Dawson function and of the scaled complementary error function to compute the moments in this limit (see \ref{ann2} for details).
\end{rema}

\section{Numerical results} \label{sec_res}
In this section, a numerical study of the scheme is presented. UGKS-M1 is compared to an asymptotic preserving modified HLL scheme for the M1 model \cite{turpault} and to the kinetic solution given by the UGKS. UGKS-M2 is also tested and compared. We chose different test cases to validate all the regimes and the convergence order. The simulation parameters are summarized in table \ref{tab:param}. 

\begin{table}[!htbp] 
     \center
    \begin{tabular}{| c || l | l | l | l | l | l | l |}
     \hline			
        ~ & $\eta$ & $\epsilon$ & $\sigma(x)$ & $f(t,0,v)$ & $f(t,1,v)$ & $\rho(0,x)$ & $u(0,x)$\\\hline
       Convergence & 1 & 1 & 1 & Periodic & Periodic & $0.5+0.25\sin{(2\pi x)}$ & 0.4\\
       Transport & 1 & 1 & 1 & 0 & $\mathbb{1}_{v<0}$ & 0 & 0 \\
       Intermediate & $10^{-1}$ & $10^{-1}$ & 1 & 0 & $\mathbb{1}_{v<0}$ & 0 & 0\\
       Diffusion & $10^{-8}$ & $10^{-8}$ & 1 & $\mathbb{1}_{v>0}$ & 0 & 0 & 0\\
     \hline  
     \end{tabular}
     \caption{Simulation parameters\label{tab:param}.}
\end{table}
 The spatial domain $\mathscr{D}=[0,1]$ is discretized with 200 points and the velocity space with 50 points (for the UGKS). Two types of boundary conditions are considered: the Dirichlet condition where the distribution function is enforced at the boundary and the periodic condition.

\paragraph{Test n°1: Relaxation of a sinusoid in a infinite domain.} 

\

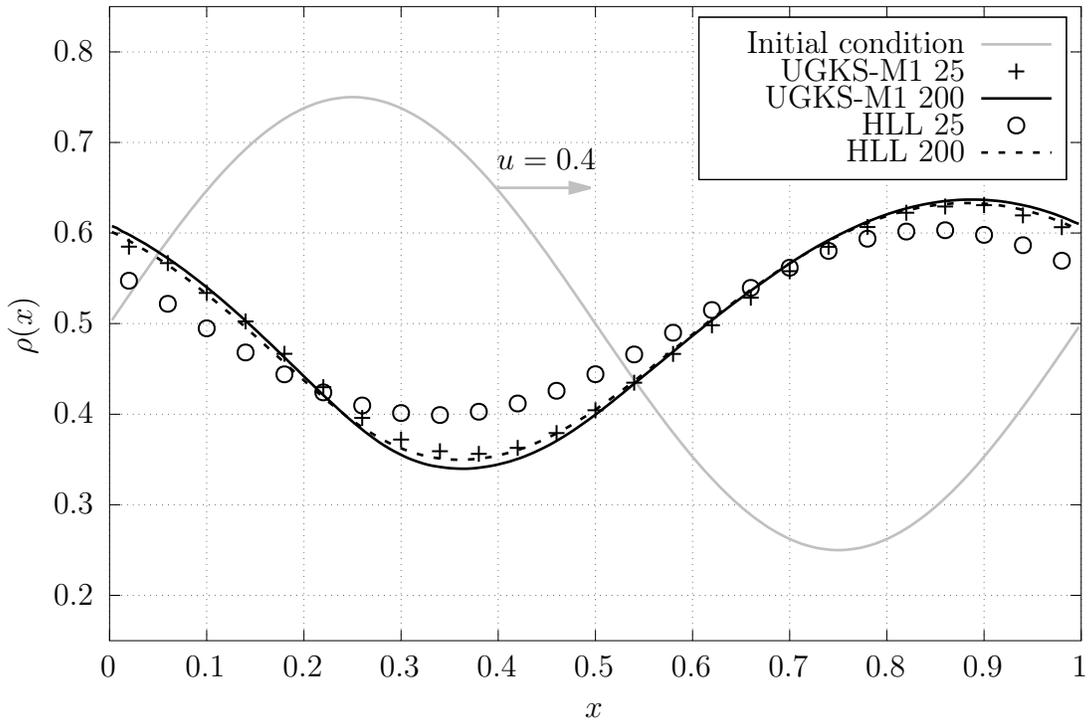
\begin{figure}[!htbp]
    \centering \input{Aperiodique}
    \caption{Test n°1: mesh convergence study for UGKS-M1 and HLL. Density in the domain at time $t=1.0$.}
    \label{fig:cinetique_conv}
\end{figure}

Firstly, a regular initial condition is considered with a sinusoidal density distribution and a uniform velocity. The periodic boundary conditions are equivalent to the transport of the sinusoid in an infinite domain. 
In figure \ref{fig:cinetique_conv}, we observe that the density is mostly advected to the right. Moreover, the amplitude of the sine wave is reduced by 15\% due to the diffusion involved by the relaxation towards the equilibrium. From a numerical point of view, we notice that UGKS-M1 solution has less diffusion than the HLL one. This phenomenon is a consequence of the choice of the waves speeds in the approximate Riemann solver. A standard choice is to use the extreme values of the Jacobian eigenvalues. However since the M1 moment closure is not analytical (for the Boltzmann entropy), we have chosen to bound those values. This choice induces numerical diffusion. 

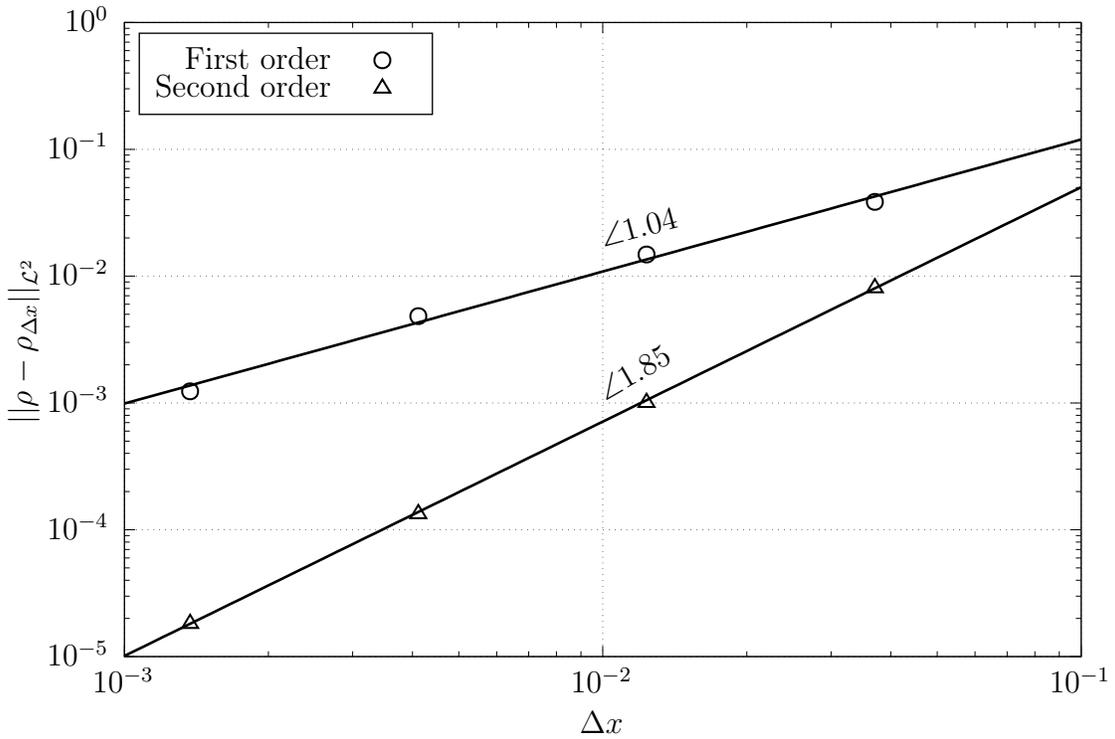
\begin{figure}[!htbp]
    \centering \input{Aconvergence}
    \caption{Test n°1: UGKS-M1 density error as a function of the step size.}
    \label{fig:conv}
\end{figure}

 For this test case and with the first and second order version of UGKS-M1, we plot in figure \ref{fig:conv} the $\mathrm{L}^2$ norm of the density error ${\Tilde\rho - \rho_{\Delta x}}$ against the step size. The reference solution $\Tilde\rho$ is computed on a grid that is small enough to assume that the error in relation to the exact solution is negligible in the analysis. The Van Leer limiter is used \cite{vanleer}. A linear regression allows to compute the convergence order of both schemes. The linear reconstruction with slope limiter leads to a significantly higher order of $1.85$ on this test case.

\paragraph{Test n°2: Transport regime with Dirichlet boundary conditions} 

\

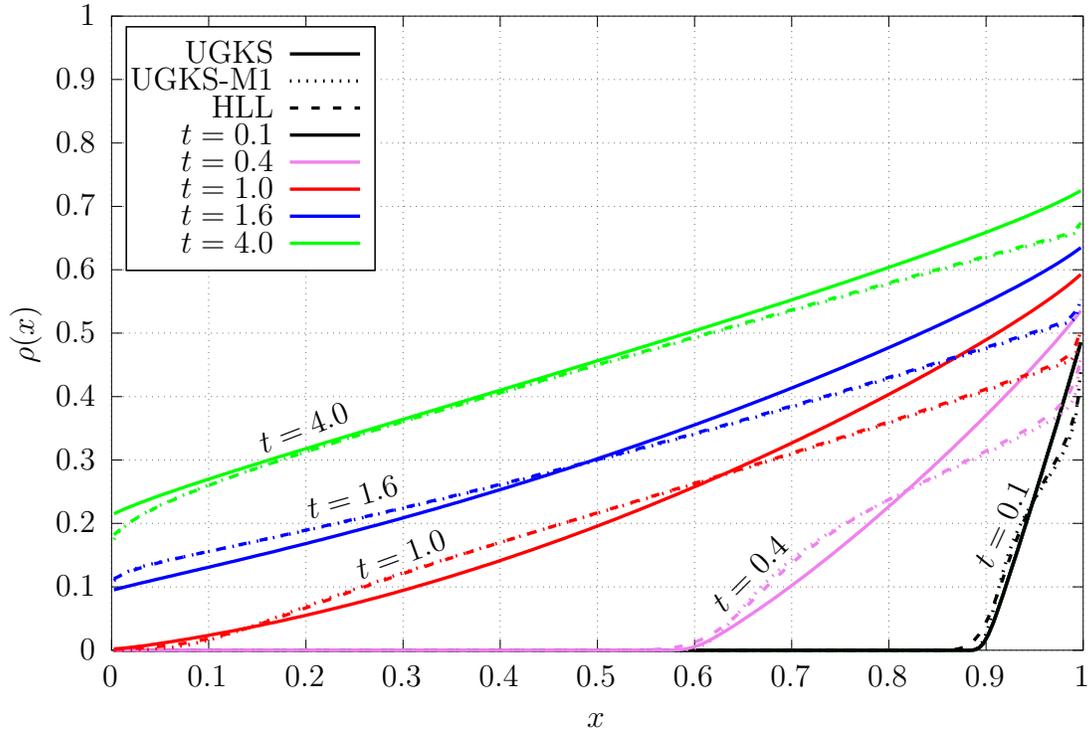
\begin{figure}[!htbp]
    \centering \input{Acinetique}
    \caption{Test n°2: transport regime. Density in the domain at different times for UGKS and M1 (with UGKS-M1 and HLL).}
    \label{fig:cinetique}
\end{figure}

In this test case, we consider a null density initial condition. On both sides of the domain, a uniform half distribution function is enforced at the kinetic level for entering particles. For the UGKS and therefore UGKS-M1, the numerical flux at the boundary is obtained by modifying the distribution function representation at the boundary by setting (for example at the left boundary):
\begin{equation}
    \label{representation_interface_bord}
f(t,x_{1/2},v)=
    \begin{cases}
    f_L(t,v) & \text{ if } v>0 \\
    \begin{array}{l}
e^{-\nu_{1/2}(t-t_n)} f(t_n,x_{1/2}-\dfrac{v}{\eta}(t-t_n),v) \\
+\nu_{1/2}\displaystyle \int_{t_n}^t  e^{-\nu_{1/2}(t-s)} \rho(s,x_{1/2} - \dfrac{v}{\eta}(t-s)) \mathrm{d}s 
\end{array}& \text{ if } v<0
    \end{cases}.
\end{equation}
Thus, the microscopic flux is:
\begin{equation}
    \phi_{1/2}= \dfrac{v}{\eta} f_L \mathbb{1}_{v>0} + (
    A_{1/2} v f_1^n + C_{1/2} v \rho_{1/2}^n
     + D_{1/2} v^2 \delta_x^L \rho_{1/2}^{n} \mathbb{1} _{v<0}
    )\mathbb{1}_{v<0},
\end{equation}
where the interface density is artificially set to $\rho_{1/2}=-\frac{\crochet{vf_L\mathbb{1}_{v>0}}}{\crochet{v\mathbb{1}_{v<0}}}$ to ensure a good asymptotic behavior \cite{mieussens2013}. For the HLL scheme, a ghost cell is used to implement Marshak boundary conditions:
\begin{equation}
\mathbf{U}_{0}^n =  \begin{pmatrix}
    \crochet{f_L \mathbb{1}_{v>0} + \hat{f}_1^n \mathbb{1}_{v<0}} \\
    \crochet{vf_L \mathbb{1}_{v>0}+ v\hat{f}_1^n \mathbb{1}_{v<0}} \\
\end{pmatrix}.  
\end{equation}
In figure \ref{fig:cinetique}, the density in the domain is represented at different times. Both M1 solutions are almost indistinguishable. Before $t=0.4$, we can notice that HLL is still slightly more diffusive than UGKS-M1 especially near the boundary and the front of the wave. The distribution function becomes isotropic over time, and the density reaches a stationary regular state. In that limit, the two computed densities tend to be identical. 

 The UGKS solution is significantly different from both M1 solutions. This is due to the fact that in the transport regime, the M1 model is highly inaccurate as compared to the underlying kinetic equation. Indeed, the distribution functions are highly out of equilibrium, thus the projections on the M1 set are inaccurate. For example at the right boundary, the density is systematically lower at all times because the projection of the half distribution function at the boundary leads to the creation of positive velocity particles. As a consequence, fewer particles enter the domain. The isotropization process alleviates this problem over time.

\paragraph{Test n°3: Intermediate regime with Dirichlet boundary conditions}

\

In figure \ref{fig:intermediaire}, we can notice that at a lower Knudsen number, both M1 solutions are much closer to the solution of the kinetic equation. In intermediate regimes the M1 model is much more relevant as the distribution functions are rapidly close to the equilibrium. UGKS-M1 is almost indistinguishable from the UGKS except at $t=0.1$ where the amplitude is 2\% lower close to the boundary. The HLL solution has again more numerical diffusion; before reaching the stationary state a significant gap can be observed in the whole domain. 

\paragraph{Test n°4: Diffusion regime with Dirichlet boundary conditions} 

\

In the diffusion regime, the solutions are identical as all schemes degenerate in the same way (as shown in figure \ref{fig:diffusion}).

\paragraph{Test n°5: UGKS-M2} 

\

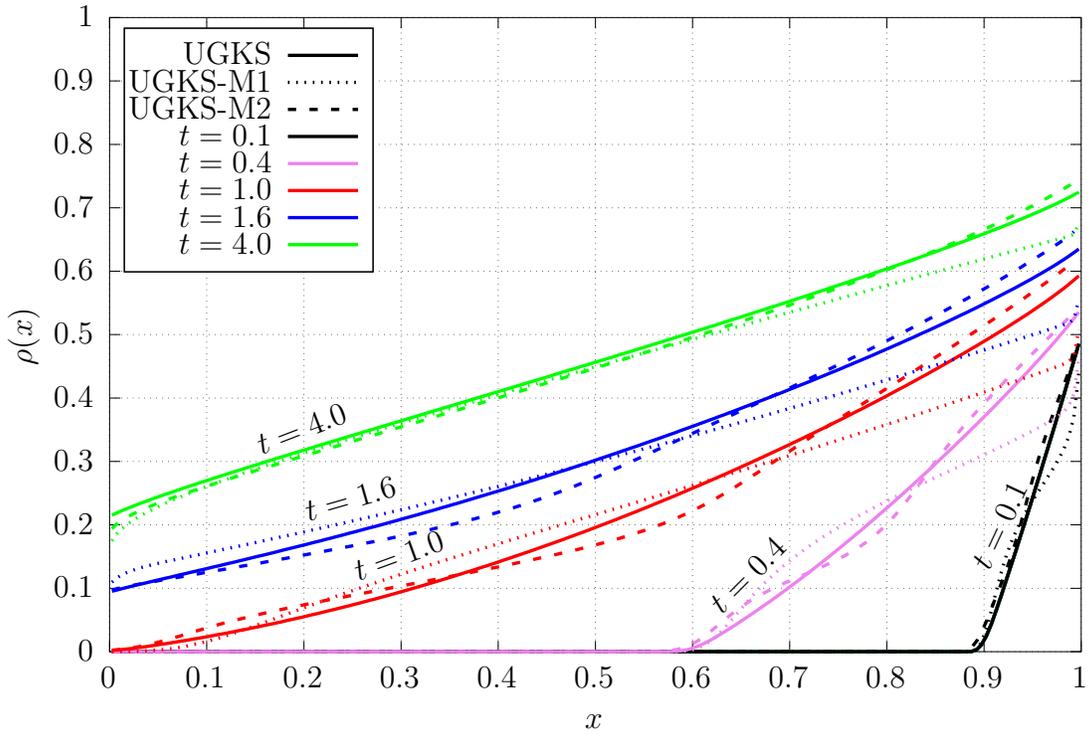
\begin{figure}[!thbp]
    \centering \input{Acinetique_bis}
    \caption{Test n°5: UGKS-M2. Density in the domain at different times for UGKS, UGKS-M1 and UGKS-M2 in the transport regime.}
    \label{fig:cinetique_bis}
\end{figure}

In test case n°3, the results emphasized that the M1 model is highly inaccurate in the transport regime. The M2 moment model allows to increase the order of the hierarchy and consequently expands the set of representable distribution functions. Thus, this model should be able to recover more kinetic effects. In figure \ref{fig:cinetique_bis}, the solution of the M2 moment model (given by UGKS-M2) is compared to the M1 solution and to the kinetic one in the transport regime. Overall, a significant improvement in the results can be noticed compared to the M1 solutions. At the boundary, the gap between the kinetic solution is much less pronounced as the enforced distribution function is better represented. In the domain, the density is also closer to the kinetic one, especially at times $t=0.1$ and $t=0.4$ when the distribution is highly out of equilibrium. Although the solution is better, the M2 model is still not sufficient to perfectly solve this kinetic case. 

\begin{figure}[!thbp]
    \centering \input{Aintermediaire_bis}
    \caption{Test n°5: UGKS-M2. Density in the domain at different times for UGKS, UGKS-M1 and UGKS-M2 in an intermediate regime.}
    \label{fig:intermediate_bis}
\end{figure}
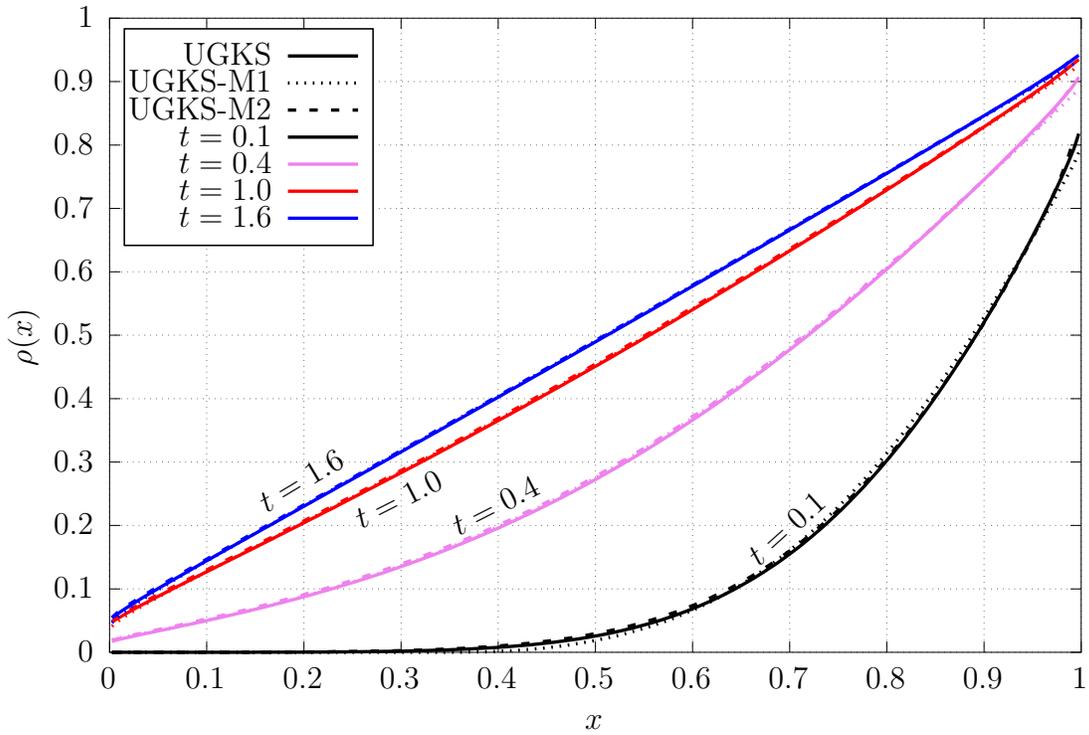

In figure \ref{fig:intermediate_bis}, UGKS-M2 is compared to UGKS-M1 and UGKS in the intermediate regime. In this test case, the difference with the kinetic solution is almost indistinguishable. At the first time $t=0.1$, it can be noted that the M2 solution is correct at the boundary unlike the M1 solution.

UGKS-M2 has also been tested and compared in the diffusion regime and gives exactly the same perfect results as M1 (see test case n°4) and are not shown here.

\section{Conclusion}
In this work, an asymptotic-preserving scheme based on the Unified Gas Kinetic Scheme has been proposed for the M1 model of linear transport. This new method consists in performing a numerical moment closure in the UGKS fluxes using the M1 distribution function. It has been demonstrated that this procedure allows to inherit the asymptotic-preserving property of the UGKS and hence to recover correct numerical fluxes in the diffusion limit.
Moreover, a second-order extension that does not compromise the AP property has been suggested. Several test cases have been chosen to validate and showcase the good behavior of the scheme in all regimes. This method has also been compared with a HLL asymptotic-preserving scheme and proved to be more accurate, especially in intermediate regimes. Furthermore, the generic nature of this method has been demonstrated with an application to the M2 moment model. Finally, this second scheme is tested and the advantages of this model as compared to the M1 one are highlighted in a kinetic case.

Despite the rather simple physical context, this article proposes a general procedure to obtain good asymptotic-preserving schemes for moment models. A work currently underway focuses on ensuring the preservation of the admissible states under CFL-like conditons in UGKS-M1. This property could be obtained by modifying the UGKS to ensure the positivity of the distribution function at kinetic level and hence ensure that the M1 variables are moments of a positive distribution function. Another relevant perspective would be to apply this procedure on other moment models based on more relevant collision kernels in higher dimension and on unstructured meshes. It would also be interesting to study non-linear collision operators with this approach. 

\bibliography{biblio}
\clearpage

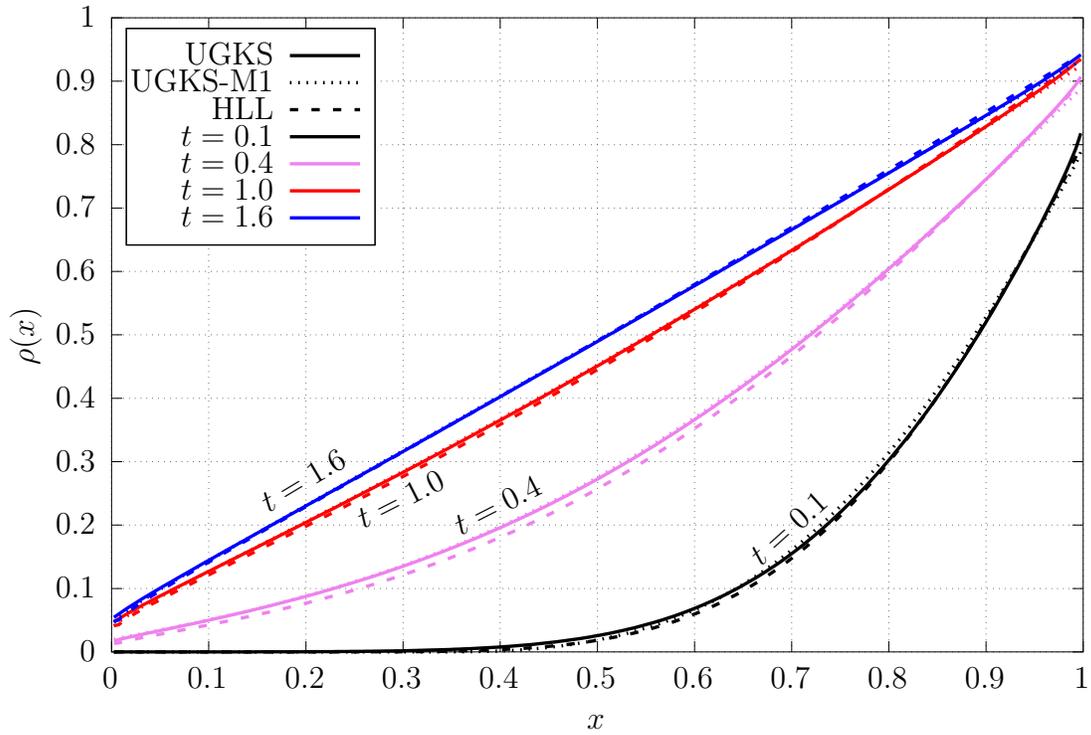
\begin{figure}[!htbp]
    \centering \input{Aintermediaire}
    \caption{Test n°3: intermediate regime. Density in the domain at different times for UGKS and M1 (with UGKS-M1 and HLL).}
    \label{fig:intermediaire}
\end{figure}
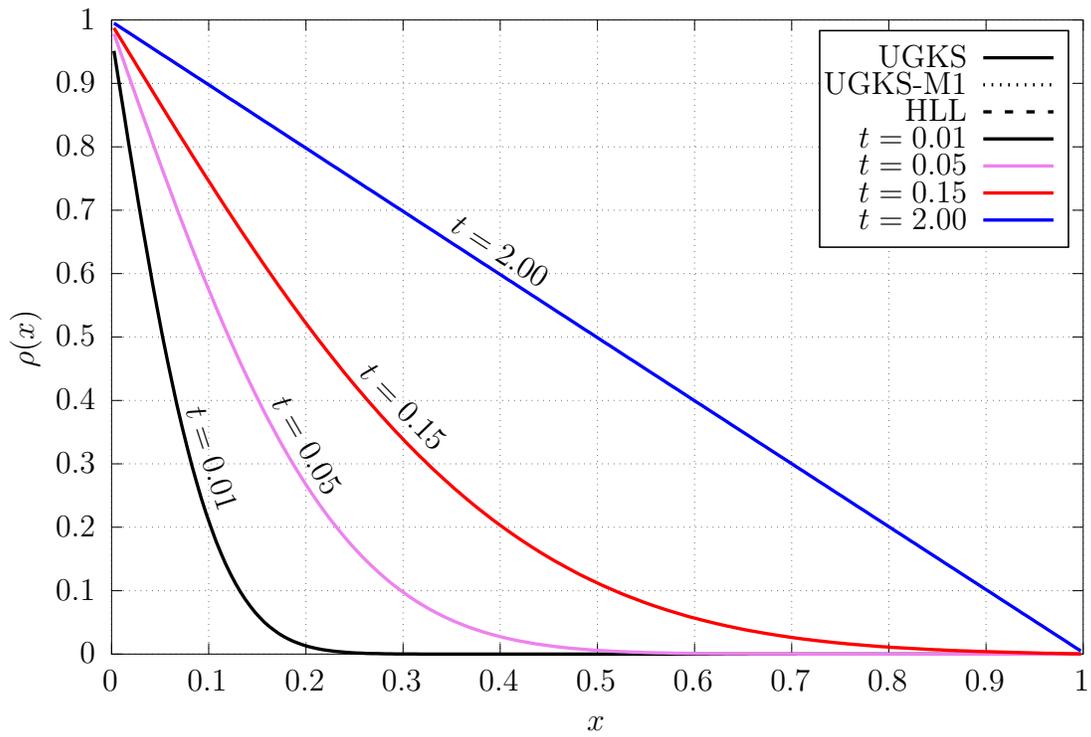
\begin{figure}[!thbp]
    \centering \input{Adiffusion}
    \caption{Test n°4: diffusion regime. Density in the domain at different times for UGKS and M1 (with UGKS-M1 and HLL).}
    \label{fig:diffusion}
\end{figure}

\clearpage

\clearpage
\appendix
\section{UGKS-M1 fluxes \label{ann1}}
In developed form, the second order numerical fluxes of UGKS-M1 are
\begin{subequations}
    \begin{align}
        \begin{split}
            {{\Phi}_{i+1/2}^\rho} &= A_{i+1/2}\left[\crochet{v\hat{f}_{i}^n\mathbb{1}_{v>0}}+\crochet{v\hat{f}_{i+1}^n\mathbb{1}_{v<0}}\right]\\
            &+A_{i+1/2} \dfrac{\Delta x}{2}\left[\mathbf{J}_\mathbf\Lambda(\mathbf{U}_i^n) \mathbf{\delta U}_i^n\cdot\begin{pmatrix} \crochet{v\hat{f}_{i}^n\mathbb{1}_{v>0}} \\ \crochet{v^2\hat{f}_{i}^n\mathbb{1}_{v>0}} \end{pmatrix} -  \mathbf{J}_\mathbf\Lambda(\mathbf{U}_{i+1}^n) \mathbf{\delta U}_{i+1}^n\cdot\begin{pmatrix} \crochet{v\hat{f}_{i+1}^n\mathbb{1}_{v<0}} \\ \crochet{v^2\hat{f}_{i+1}^n\mathbb{1}_{v<0}} \end{pmatrix} \right] \\
            &+B_{i+1/2}\left[\mathbf{J}_\mathbf\Lambda(\mathbf{U}_i^n) \mathbf{\delta U}_i^n\cdot\begin{pmatrix} \crochet{v^2\hat{f}_{i}^n\mathbb{1}_{v>0}} \\ \crochet{v^3\hat{f}_{i}^n\mathbb{1}_{v>0}} \end{pmatrix} +  \mathbf{J}_\mathbf\Lambda(\mathbf{U}_{i+1}^n) \mathbf{\delta U}_{i+1}^n\cdot\begin{pmatrix} \crochet{v^2\hat{f}_{i+1}^n\mathbb{1}_{v<0}} \\ \crochet{v^3\hat{f}_{i+1}^n\mathbb{1}_{v<0}} \end{pmatrix} \right] \\
            &+ \dfrac{D_{i+1/2}}{3\Delta x} (\rho_{i+1}^{n}-\rho_i^n), \\
        \end{split} \\
        \begin{split}
            {{\Phi}_{i+1/2}^j} &=A_{i+1/2}\left[\crochet{v^2\hat{f}_{i}^n\mathbb{1}_{v>0}}+\crochet{v^2\hat{f}_{i+1}^n\mathbb{1}_{v<0}}\right]\\
            &+A_{i+1/2} \dfrac{\Delta x}{2}\left[\mathbf{J}_\mathbf\Lambda(\mathbf{U}_i^n) \mathbf{\delta U}_i^n\cdot\begin{pmatrix} \crochet{v^2\hat{f}_{i}^n\mathbb{1}_{v>0}} \\ \crochet{v^3\hat{f}_{i}^n\mathbb{1}_{v>0}} \end{pmatrix} -  \mathbf{J}_\mathbf\Lambda(\mathbf{U}_{i+1}^n) \mathbf{\delta U}_{i+1}^n\cdot\begin{pmatrix} \crochet{v^2\hat{f}_{i+1}^n\mathbb{1}_{v<0}} \\ \crochet{v^3\hat{f}_{i+1}^n\mathbb{1}_{v<0}} \end{pmatrix} \right] \\
            &+B_{i+1/2}\left[\mathbf{J}_\mathbf\Lambda(\mathbf{U}_i^n) \mathbf{\delta U}_i^n\cdot\begin{pmatrix} \crochet{v^3\hat{f}_{i}^n\mathbb{1}_{v>0}} \\ \crochet{v^4\hat{f}_{i}^n\mathbb{1}_{v>0}} \end{pmatrix}+  \mathbf{J}_\mathbf\Lambda(\mathbf{U}_{i+1}^n) \mathbf{\delta U}_{i+1}^n\cdot\begin{pmatrix} \crochet{v^3\hat{f}_{i+1}^n\mathbb{1}_{v<0}} \\ \crochet{v^4\hat{f}_{i+1}^n\mathbb{1}_{v<0}} \end{pmatrix} \right] \\
            &+ \dfrac{C_{i+1/2}}{3}\left[\crochet{\hat{f}_{i}^n\mathbb{1}_{v>0}}+\crochet{\hat{f}_{i+1}^n\mathbb{1}_{v<0}}\right].\\
        \end{split}
    \end{align}
\end{subequations}
The half-moments of the M1 distribution function are expressed in developed form to reduce the accumulation of round-off error at low $\beta$: 
\begin{subequations}
\begin{align*}
    \crochet{\hat{f}_i^n \mathbb{1}_{v>0}}&=\dfrac{e^{\alpha_i}}{2\beta_i}\left(e^{\beta_i}-1\right), \\
    \crochet{v\hat{f}_i^n \mathbb{1}_{v>0}}&=\dfrac{e^{\alpha_i}}{2\beta_i}\left(e^{\beta_i}\left(1-\dfrac{1}{\beta_i}\right)+\dfrac{1}{\beta_i}\right),\\
    \crochet{v^2\hat{f}_i^n \mathbb{1}_{v>0}}&=\dfrac{e^{\alpha_i}}{2\beta_i}\left(e^{\beta_i}\left(1-\dfrac{2}{\beta_i}+\dfrac{2}{\beta_i^2}\right)-\dfrac{2}{\beta_i^2}\right), \\
    \crochet{v^3\hat{f}_i^n \mathbb{1}_{v>0}}&=\dfrac{e^{\alpha_i}}{2\beta_i}\left(e^{\beta_i}\left(1-\dfrac{3}{\beta_i}+\dfrac{6}{\beta_i^2}-\dfrac{6}{\beta_i^3}\right)+\dfrac{6}{\beta_i^3}\right), \\
    \crochet{v^4\hat{f}_i^n \mathbb{1}_{v>0}}&=\dfrac{e^{\alpha_i}}{2\beta_i}\left(e^{\beta_i}\left(1-\dfrac{4}{\beta_i}+\dfrac{12}{\beta_i^2}-\dfrac{24}{\beta_i^3}+\dfrac{24}{\beta_i^4}\right)-\dfrac{24}{\beta_i^4}\right),
\end{align*}
\end{subequations}
\begin{subequations}
\begin{align*}
    \crochet{\hat{f}_{i}^n \mathbb{1}_{v<0}}&=\dfrac{-e^{\alpha_{i}}}{2\beta_{i}}\left(e^{-\beta_{i}}-1\right), \\
    \crochet{v\hat{f}_{i}^n \mathbb{1}_{v<0}}&=\dfrac{e^{\alpha_{i}}}{2\beta_{i}}\left(e^{-\beta_{i}}\left(1+\dfrac{1}{\beta_i}\right)-\dfrac{1}{\beta_i}\right), \\
    \crochet{v^2\hat{f}_i^n \mathbb{1}_{v<0}}&=\dfrac{-e^{\alpha_i}}{2\beta_i}\left(e^{-\beta_i}\left(1+\dfrac{2}{\beta_i}+\dfrac{2}{\beta_i^2}\right)-\dfrac{2}{\beta_i^2}\right), \\
    \crochet{v^3\hat{f}_i^n \mathbb{1}_{v<0}}&=\dfrac{e^{\alpha_i}}{2\beta_i}\left(e^{-\beta_i}\left(1+\dfrac{3}{\beta_i}+\dfrac{6}{\beta_i^2}+\dfrac{6}{\beta_i^3}\right)-\dfrac{6}{\beta_i^3}\right), \\
    \crochet{v^4\hat{f}_i^n \mathbb{1}_{v<0}}&=\dfrac{-e^{\alpha_i}}{2\beta_i}\left(e^{-\beta_i}\left(1+\dfrac{4}{\beta_i}+\dfrac{12}{\beta_i^2}+\dfrac{24}{\beta_i^3}+\dfrac{24}{\beta_i^4}\right)-\dfrac{24}{\beta_i^4}\right). 
\end{align*}
\end{subequations}

\section{Computation of the M2 density \label{ann2}}
As $\frac{\beta}{2\sqrt\gamma}$ tends to infinity, the computation of the density (\ref{p_m2}) and the half densities (\ref{half_densities_m2}) of the M2 distribution function becomes numerically stiff. For example, in the case $\gamma<0$ and $\beta>0$, the quantities inside the density expression are not bounded. We found that a robust computation in this limit is to use the asymptotic expansions of the Dawson function $\Dp$ and of the scaled complementary error function $\Dm$:
\begin{subequations}
\begin{align}
\text{For large }x,~\Dp(x)&=\dfrac{1}{2x}\left(1+\sum_{j=1}^n \dfrac{(2j-1)!!}{2^jx^{2j}}\right) = \dfrac{1}{2x}\left(1+ \dfrac{1}{2x^2}+\dfrac{3}{4x^4}+\cdots\right), \\
\text{For large }x,~\Dm(x)&=\dfrac{1}{2x}\left(1+\sum_{j=1}^n (-1)^j \dfrac{(2j-1)!!}{2^jx^{2j}}\right) + \sqrt\pi e^{x^2}\mathbb{1}_{x\leq0}.
\end{align}
\end{subequations}
where $n$ is the number of terms in the series. This number should be greater than $3$ to recover every M2 moments but not too large as the series expansions of both functions are non convergent. Those same expansions can be used to formally recover the M1 moments from the M2 ones.
\end{document}

%% file: Aperiodique.tex
\begingroup
  \makeatletter
  \providecommand\color[2][]{%
    \GenericError{(gnuplot) \space\space\space\@spaces}{%
      Package color not loaded in conjunction with
      terminal option `colourtext'%
    }{See the gnuplot documentation for explanation.%
    }{Either use 'blacktext' in gnuplot or load the package
      color.sty in LaTeX.}%
    \renewcommand\color[2][]{}%
  }%
  \providecommand\includegraphics[2][]{%
    \GenericError{(gnuplot) \space\space\space\@spaces}{%
      Package graphicx or graphics not loaded%
    }{See the gnuplot documentation for explanation.%
    }{The gnuplot epslatex terminal needs graphicx.sty or graphics.sty.}%
    \renewcommand\includegraphics[2][]{}%
  }%
  \providecommand\rotatebox[2]{#2}%
  \@ifundefined{ifGPcolor}{%
    \newif\ifGPcolor
    \GPcolortrue
  }{}%
  \@ifundefined{ifGPblacktext}{%
    \newif\ifGPblacktext
    \GPblacktexttrue
  }{}%
  \let\gplgaddtomacro\g@addto@macro
  \gdef\gplbacktext{}%
  \gdef\gplfronttext{}%
  \makeatother
  \ifGPblacktext
    \def\colorrgb#1{}%
    \def\colorgray#1{}%
  \else
    \ifGPcolor
      \def\colorrgb#1{\color[rgb]{#1}}%
      \def\colorgray#1{\color[gray]{#1}}%
      \expandafter\def\csname LTw\endcsname{\color{white}}%
      \expandafter\def\csname LTb\endcsname{\color{black}}%
      \expandafter\def\csname LTa\endcsname{\color{black}}%
      \expandafter\def\csname LT0\endcsname{\color[rgb]{1,0,0}}%
      \expandafter\def\csname LT1\endcsname{\color[rgb]{0,1,0}}%
      \expandafter\def\csname LT2\endcsname{\color[rgb]{0,0,1}}%
      \expandafter\def\csname LT3\endcsname{\color[rgb]{1,0,1}}%
      \expandafter\def\csname LT4\endcsname{\color[rgb]{0,1,1}}%
      \expandafter\def\csname LT5\endcsname{\color[rgb]{1,1,0}}%
      \expandafter\def\csname LT6\endcsname{\color[rgb]{0,0,0}}%
      \expandafter\def\csname LT7\endcsname{\color[rgb]{1,0.3,0}}%
      \expandafter\def\csname LT8\endcsname{\color[rgb]{0.5,0.5,0.5}}%
    \else
      \def\colorrgb#1{\color{black}}%
      \def\colorgray#1{\color[gray]{#1}}%
      \expandafter\def\csname LTw\endcsname{\color{white}}%
      \expandafter\def\csname LTb\endcsname{\color{black}}%
      \expandafter\def\csname LTa\endcsname{\color{black}}%
      \expandafter\def\csname LT0\endcsname{\color{black}}%
      \expandafter\def\csname LT1\endcsname{\color{black}}%
      \expandafter\def\csname LT2\endcsname{\color{black}}%
      \expandafter\def\csname LT3\endcsname{\color{black}}%
      \expandafter\def\csname LT4\endcsname{\color{black}}%
      \expandafter\def\csname LT5\endcsname{\color{black}}%
      \expandafter\def\csname LT6\endcsname{\color{black}}%
      \expandafter\def\csname LT7\endcsname{\color{black}}%
      \expandafter\def\csname LT8\endcsname{\color{black}}%
    \fi
  \fi
    \setlength{\unitlength}{0.0500bp}%
    \ifx\gptboxheight\undefined%
      \newlength{\gptboxheight}%
      \newlength{\gptboxwidth}%
      \newsavebox{\gptboxtext}%
    \fi%
    \setlength{\fboxrule}{0.5pt}%
    \setlength{\fboxsep}{1pt}%
\begin{picture}(8480.00,5640.00)%
    \gplgaddtomacro\gplbacktext{%
      \csname LTb\endcsname
      \put(708,994){\makebox(0,0)[r]{\strut{}$0.2$}}%
      \csname LTb\endcsname
      \put(708,1677){\makebox(0,0)[r]{\strut{}$0.3$}}%
      \csname LTb\endcsname
      \put(708,2360){\makebox(0,0)[r]{\strut{}$0.4$}}%
      \csname LTb\endcsname
      \put(708,3044){\makebox(0,0)[r]{\strut{}$0.5$}}%
      \csname LTb\endcsname
      \put(708,3727){\makebox(0,0)[r]{\strut{}$0.6$}}%
      \csname LTb\endcsname
      \put(708,4410){\makebox(0,0)[r]{\strut{}$0.7$}}%
      \csname LTb\endcsname
      \put(708,5093){\makebox(0,0)[r]{\strut{}$0.8$}}%
      \csname LTb\endcsname
      \put(820,448){\makebox(0,0){\strut{}$0$}}%
      \csname LTb\endcsname
      \put(1552,448){\makebox(0,0){\strut{}$0.1$}}%
      \csname LTb\endcsname
      \put(2285,448){\makebox(0,0){\strut{}$0.2$}}%
      \csname LTb\endcsname
      \put(3017,448){\makebox(0,0){\strut{}$0.3$}}%
      \csname LTb\endcsname
      \put(3749,448){\makebox(0,0){\strut{}$0.4$}}%
      \csname LTb\endcsname
      \put(4482,448){\makebox(0,0){\strut{}$0.5$}}%
      \csname LTb\endcsname
      \put(5214,448){\makebox(0,0){\strut{}$0.6$}}%
      \csname LTb\endcsname
      \put(5946,448){\makebox(0,0){\strut{}$0.7$}}%
      \csname LTb\endcsname
      \put(6678,448){\makebox(0,0){\strut{}$0.8$}}%
      \csname LTb\endcsname
      \put(7411,448){\makebox(0,0){\strut{}$0.9$}}%
      \csname LTb\endcsname
      \put(8143,448){\makebox(0,0){\strut{}$1$}}%
      \csname LTb\endcsname
      \put(3749,4273){\makebox(0,0)[l]{\strut{}$u=0.4$}}%
    }%
    \gplgaddtomacro\gplfronttext{%
      \csname LTb\endcsname
      \put(186,3043){\rotatebox{-270}{\makebox(0,0){\strut{}$\rho(x)$}}}%
      \csname LTb\endcsname
      \put(4481,142){\makebox(0,0){\strut{}$x$}}%
      \csname LTb\endcsname
      \put(7278,5150){\makebox(0,0)[r]{\strut{}Initial condition}}%
      \csname LTb\endcsname
      \put(7278,4946){\makebox(0,0)[r]{\strut{}UGKS-M1 25}}%
      \csname LTb\endcsname
      \put(7278,4742){\makebox(0,0)[r]{\strut{}UGKS-M1 200}}%
      \csname LTb\endcsname
      \put(7278,4538){\makebox(0,0)[r]{\strut{}HLL 25}}%
      \csname LTb\endcsname
      \put(7278,4334){\makebox(0,0)[r]{\strut{}HLL 200}}%
    }%
    \gplbacktext
    \put(0,0){\includegraphics{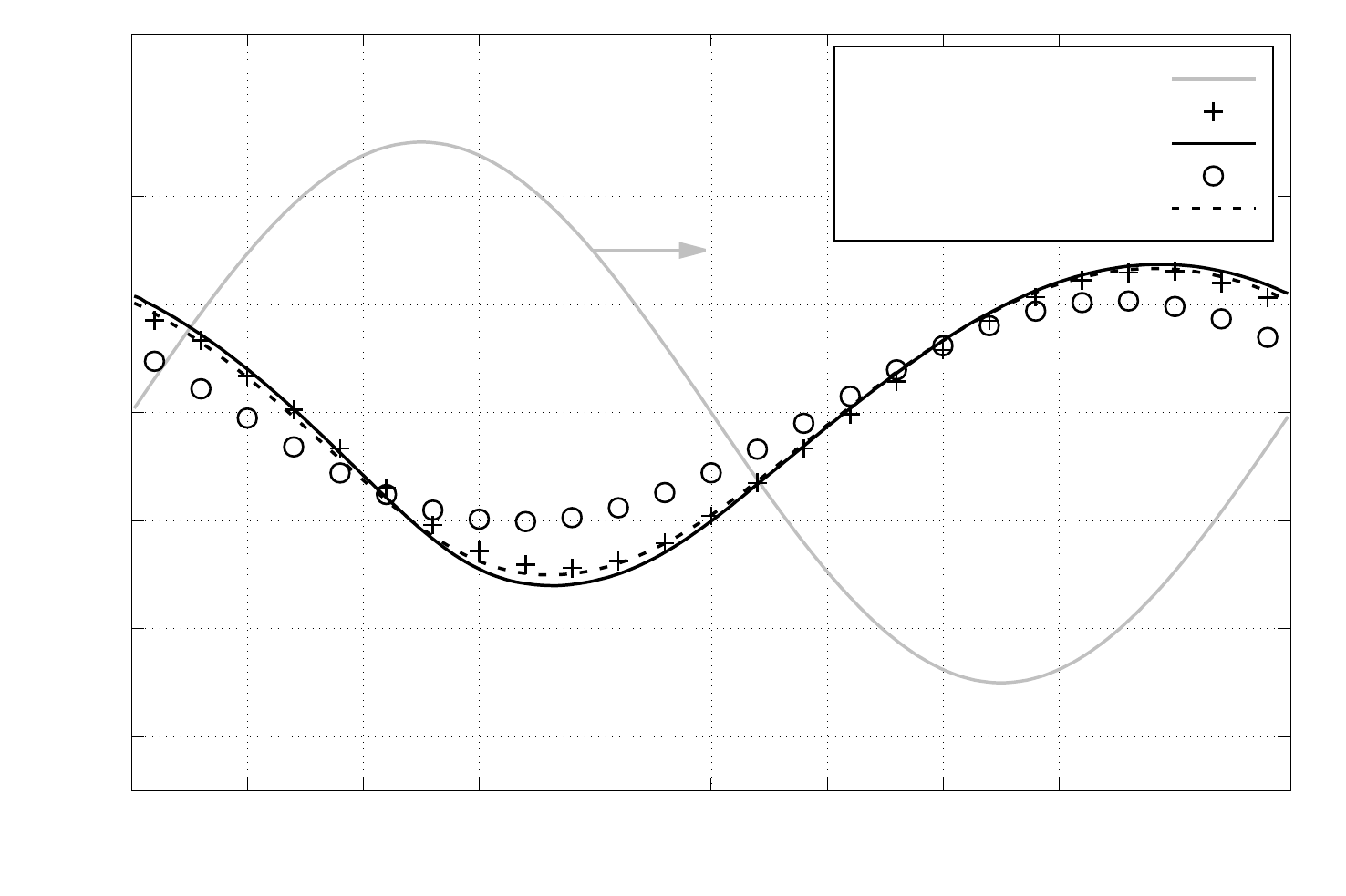}}%
    \gplfronttext
  \end{picture}%
\endgroup

%% file: Aconvergence.tex
\begingroup
  \makeatletter
  \providecommand\color[2][]{%
    \GenericError{(gnuplot) \space\space\space\@spaces}{%
      Package color not loaded in conjunction with
      terminal option `colourtext'%
    }{See the gnuplot documentation for explanation.%
    }{Either use 'blacktext' in gnuplot or load the package
      color.sty in LaTeX.}%
    \renewcommand\color[2][]{}%
  }%
  \providecommand\includegraphics[2][]{%
    \GenericError{(gnuplot) \space\space\space\@spaces}{%
      Package graphicx or graphics not loaded%
    }{See the gnuplot documentation for explanation.%
    }{The gnuplot epslatex terminal needs graphicx.sty or graphics.sty.}%
    \renewcommand\includegraphics[2][]{}%
  }%
  \providecommand\rotatebox[2]{#2}%
  \@ifundefined{ifGPcolor}{%
    \newif\ifGPcolor
    \GPcolortrue
  }{}%
  \@ifundefined{ifGPblacktext}{%
    \newif\ifGPblacktext
    \GPblacktexttrue
  }{}%
  \let\gplgaddtomacro\g@addto@macro
  \gdef\gplbacktext{}%
  \gdef\gplfronttext{}%
  \makeatother
  \ifGPblacktext
    \def\colorrgb#1{}%
    \def\colorgray#1{}%
  \else
    \ifGPcolor
      \def\colorrgb#1{\color[rgb]{#1}}%
      \def\colorgray#1{\color[gray]{#1}}%
      \expandafter\def\csname LTw\endcsname{\color{white}}%
      \expandafter\def\csname LTb\endcsname{\color{black}}%
      \expandafter\def\csname LTa\endcsname{\color{black}}%
      \expandafter\def\csname LT0\endcsname{\color[rgb]{1,0,0}}%
      \expandafter\def\csname LT1\endcsname{\color[rgb]{0,1,0}}%
      \expandafter\def\csname LT2\endcsname{\color[rgb]{0,0,1}}%
      \expandafter\def\csname LT3\endcsname{\color[rgb]{1,0,1}}%
      \expandafter\def\csname LT4\endcsname{\color[rgb]{0,1,1}}%
      \expandafter\def\csname LT5\endcsname{\color[rgb]{1,1,0}}%
      \expandafter\def\csname LT6\endcsname{\color[rgb]{0,0,0}}%
      \expandafter\def\csname LT7\endcsname{\color[rgb]{1,0.3,0}}%
      \expandafter\def\csname LT8\endcsname{\color[rgb]{0.5,0.5,0.5}}%
    \else
      \def\colorrgb#1{\color{black}}%
      \def\colorgray#1{\color[gray]{#1}}%
      \expandafter\def\csname LTw\endcsname{\color{white}}%
      \expandafter\def\csname LTb\endcsname{\color{black}}%
      \expandafter\def\csname LTa\endcsname{\color{black}}%
      \expandafter\def\csname LT0\endcsname{\color{black}}%
      \expandafter\def\csname LT1\endcsname{\color{black}}%
      \expandafter\def\csname LT2\endcsname{\color{black}}%
      \expandafter\def\csname LT3\endcsname{\color{black}}%
      \expandafter\def\csname LT4\endcsname{\color{black}}%
      \expandafter\def\csname LT5\endcsname{\color{black}}%
      \expandafter\def\csname LT6\endcsname{\color{black}}%
      \expandafter\def\csname LT7\endcsname{\color{black}}%
      \expandafter\def\csname LT8\endcsname{\color{black}}%
    \fi
  \fi
    \setlength{\unitlength}{0.0500bp}%
    \ifx\gptboxheight\undefined%
      \newlength{\gptboxheight}%
      \newlength{\gptboxwidth}%
      \newsavebox{\gptboxtext}%
    \fi%
    \setlength{\fboxrule}{0.5pt}%
    \setlength{\fboxsep}{1pt}%
\begin{picture}(8480.00,5640.00)%
    \gplgaddtomacro\gplbacktext{%
      \csname LTb\endcsname
      \put(820,652){\makebox(0,0)[r]{\strut{}$10^{-5}$}}%
      \csname LTb\endcsname
      \put(820,1609){\makebox(0,0)[r]{\strut{}$10^{-4}$}}%
      \csname LTb\endcsname
      \put(820,2565){\makebox(0,0)[r]{\strut{}$10^{-3}$}}%
      \csname LTb\endcsname
      \put(820,3522){\makebox(0,0)[r]{\strut{}$10^{-2}$}}%
      \csname LTb\endcsname
      \put(820,4478){\makebox(0,0)[r]{\strut{}$10^{-1}$}}%
      \csname LTb\endcsname
      \put(820,5435){\makebox(0,0)[r]{\strut{}$10^{0}$}}%
      \csname LTb\endcsname
      \put(932,448){\makebox(0,0){\strut{}$10^{-3}$}}%
      \csname LTb\endcsname
      \put(4537,448){\makebox(0,0){\strut{}$10^{-2}$}}%
      \csname LTb\endcsname
      \put(8143,448){\makebox(0,0){\strut{}$10^{-1}$}}%
    }%
    \gplgaddtomacro\gplfronttext{%
      \csname LTb\endcsname
      \put(186,3043){\rotatebox{-270}{\makebox(0,0){\strut{}${||\rho - \rho_{\Delta x}||_{\mathcal{L}^2}}$}}}%
      \csname LTb\endcsname
      \put(4537,142){\makebox(0,0){\strut{}${\Delta x}$}}%
      \csname LTb\endcsname
      \put(2500,5150){\makebox(0,0)[r]{\strut{}First order}}%
      \csname LTb\endcsname
      \put(2500,4946){\makebox(0,0)[r]{\strut{}Second order}}%
      \csname LTb\endcsname
      \put(4538,2641){\rotatebox{30}{\makebox(0,0)[l]{\strut{}$\angle 1.85$}}}%
      \csname LTb\endcsname
      \put(4538,3810){\rotatebox{15}{\makebox(0,0)[l]{\strut{}$\angle 1.04$}}}%
    }%
    \gplbacktext
    \put(0,0){\includegraphics{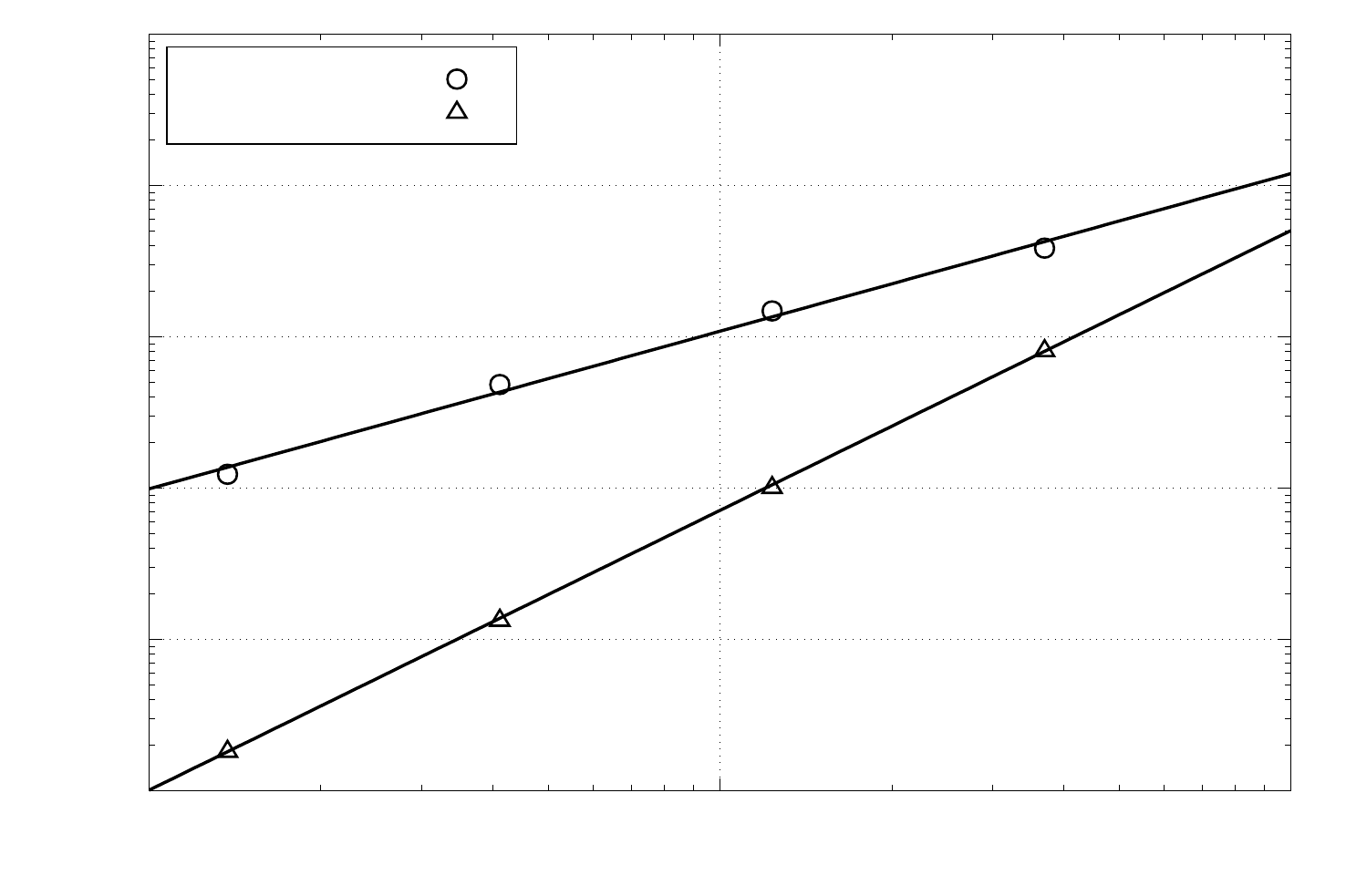}}%
    \gplfronttext
  \end{picture}%
\endgroup

%% file: Acinetique.tex
\begingroup
  \makeatletter
  \providecommand\color[2][]{%
    \GenericError{(gnuplot) \space\space\space\@spaces}{%
      Package color not loaded in conjunction with
      terminal option `colourtext'%
    }{See the gnuplot documentation for explanation.%
    }{Either use 'blacktext' in gnuplot or load the package
      color.sty in LaTeX.}%
    \renewcommand\color[2][]{}%
  }%
  \providecommand\includegraphics[2][]{%
    \GenericError{(gnuplot) \space\space\space\@spaces}{%
      Package graphicx or graphics not loaded%
    }{See the gnuplot documentation for explanation.%
    }{The gnuplot epslatex terminal needs graphicx.sty or graphics.sty.}%
    \renewcommand\includegraphics[2][]{}%
  }%
  \providecommand\rotatebox[2]{#2}%
  \@ifundefined{ifGPcolor}{%
    \newif\ifGPcolor
    \GPcolortrue
  }{}%
  \@ifundefined{ifGPblacktext}{%
    \newif\ifGPblacktext
    \GPblacktexttrue
  }{}%
  \let\gplgaddtomacro\g@addto@macro
  \gdef\gplbacktext{}%
  \gdef\gplfronttext{}%
  \makeatother
  \ifGPblacktext
    \def\colorrgb#1{}%
    \def\colorgray#1{}%
  \else
    \ifGPcolor
      \def\colorrgb#1{\color[rgb]{#1}}%
      \def\colorgray#1{\color[gray]{#1}}%
      \expandafter\def\csname LTw\endcsname{\color{white}}%
      \expandafter\def\csname LTb\endcsname{\color{black}}%
      \expandafter\def\csname LTa\endcsname{\color{black}}%
      \expandafter\def\csname LT0\endcsname{\color[rgb]{1,0,0}}%
      \expandafter\def\csname LT1\endcsname{\color[rgb]{0,1,0}}%
      \expandafter\def\csname LT2\endcsname{\color[rgb]{0,0,1}}%
      \expandafter\def\csname LT3\endcsname{\color[rgb]{1,0,1}}%
      \expandafter\def\csname LT4\endcsname{\color[rgb]{0,1,1}}%
      \expandafter\def\csname LT5\endcsname{\color[rgb]{1,1,0}}%
      \expandafter\def\csname LT6\endcsname{\color[rgb]{0,0,0}}%
      \expandafter\def\csname LT7\endcsname{\color[rgb]{1,0.3,0}}%
      \expandafter\def\csname LT8\endcsname{\color[rgb]{0.5,0.5,0.5}}%
    \else
      \def\colorrgb#1{\color{black}}%
      \def\colorgray#1{\color[gray]{#1}}%
      \expandafter\def\csname LTw\endcsname{\color{white}}%
      \expandafter\def\csname LTb\endcsname{\color{black}}%
      \expandafter\def\csname LTa\endcsname{\color{black}}%
      \expandafter\def\csname LT0\endcsname{\color{black}}%
      \expandafter\def\csname LT1\endcsname{\color{black}}%
      \expandafter\def\csname LT2\endcsname{\color{black}}%
      \expandafter\def\csname LT3\endcsname{\color{black}}%
      \expandafter\def\csname LT4\endcsname{\color{black}}%
      \expandafter\def\csname LT5\endcsname{\color{black}}%
      \expandafter\def\csname LT6\endcsname{\color{black}}%
      \expandafter\def\csname LT7\endcsname{\color{black}}%
      \expandafter\def\csname LT8\endcsname{\color{black}}%
    \fi
  \fi
    \setlength{\unitlength}{0.0500bp}%
    \ifx\gptboxheight\undefined%
      \newlength{\gptboxheight}%
      \newlength{\gptboxwidth}%
      \newsavebox{\gptboxtext}%
    \fi%
    \setlength{\fboxrule}{0.5pt}%
    \setlength{\fboxsep}{1pt}%
\begin{picture}(8480.00,5640.00)%
    \gplgaddtomacro\gplbacktext{%
      \csname LTb\endcsname
      \put(708,652){\makebox(0,0)[r]{\strut{}$0$}}%
      \csname LTb\endcsname
      \put(708,1130){\makebox(0,0)[r]{\strut{}$0.1$}}%
      \csname LTb\endcsname
      \put(708,1609){\makebox(0,0)[r]{\strut{}$0.2$}}%
      \csname LTb\endcsname
      \put(708,2087){\makebox(0,0)[r]{\strut{}$0.3$}}%
      \csname LTb\endcsname
      \put(708,2565){\makebox(0,0)[r]{\strut{}$0.4$}}%
      \csname LTb\endcsname
      \put(708,3044){\makebox(0,0)[r]{\strut{}$0.5$}}%
      \csname LTb\endcsname
      \put(708,3522){\makebox(0,0)[r]{\strut{}$0.6$}}%
      \csname LTb\endcsname
      \put(708,4000){\makebox(0,0)[r]{\strut{}$0.7$}}%
      \csname LTb\endcsname
      \put(708,4478){\makebox(0,0)[r]{\strut{}$0.8$}}%
      \csname LTb\endcsname
      \put(708,4957){\makebox(0,0)[r]{\strut{}$0.9$}}%
      \csname LTb\endcsname
      \put(708,5435){\makebox(0,0)[r]{\strut{}$1$}}%
      \csname LTb\endcsname
      \put(820,448){\makebox(0,0){\strut{}$0$}}%
      \csname LTb\endcsname
      \put(1552,448){\makebox(0,0){\strut{}$0.1$}}%
      \csname LTb\endcsname
      \put(2285,448){\makebox(0,0){\strut{}$0.2$}}%
      \csname LTb\endcsname
      \put(3017,448){\makebox(0,0){\strut{}$0.3$}}%
      \csname LTb\endcsname
      \put(3749,448){\makebox(0,0){\strut{}$0.4$}}%
      \csname LTb\endcsname
      \put(4482,448){\makebox(0,0){\strut{}$0.5$}}%
      \csname LTb\endcsname
      \put(5214,448){\makebox(0,0){\strut{}$0.6$}}%
      \csname LTb\endcsname
      \put(5946,448){\makebox(0,0){\strut{}$0.7$}}%
      \csname LTb\endcsname
      \put(6678,448){\makebox(0,0){\strut{}$0.8$}}%
      \csname LTb\endcsname
      \put(7411,448){\makebox(0,0){\strut{}$0.9$}}%
      \csname LTb\endcsname
      \put(8143,448){\makebox(0,0){\strut{}$1$}}%
    }%
    \gplgaddtomacro\gplfronttext{%
      \csname LTb\endcsname
      \put(186,3043){\rotatebox{-270}{\makebox(0,0){\strut{}$\rho(x)$}}}%
      \csname LTb\endcsname
      \put(4481,142){\makebox(0,0){\strut{}$x$}}%
      \csname LTb\endcsname
      \put(2052,5150){\makebox(0,0)[r]{\strut{}UGKS}}%
      \csname LTb\endcsname
      \put(2052,4946){\makebox(0,0)[r]{\strut{}  UGKS-M1}}%
      \csname LTb\endcsname
      \put(2052,4742){\makebox(0,0)[r]{\strut{}HLL}}%
      \csname LTb\endcsname
      \put(2052,4538){\makebox(0,0)[r]{\strut{}  $t=0.1$}}%
      \csname LTb\endcsname
      \put(2052,4334){\makebox(0,0)[r]{\strut{}  $t=0.4$}}%
      \csname LTb\endcsname
      \put(2052,4130){\makebox(0,0)[r]{\strut{}  $t=1.0$}}%
      \csname LTb\endcsname
      \put(2052,3926){\makebox(0,0)[r]{\strut{}  $t=1.6$}}%
      \csname LTb\endcsname
      \put(2052,3722){\makebox(0,0)[r]{\strut{}  $t=4.0$}}%
      \csname LTb\endcsname
      \put(7557,1609){\rotatebox{65}{\makebox(0,0){\strut{}$t=0.1$}}}%
      \csname LTb\endcsname
      \put(5653,1226){\rotatebox{45}{\makebox(0,0){\strut{}$t=0.4$}}}%
      \csname LTb\endcsname
      \put(3017,1369){\rotatebox{20}{\makebox(0,0){\strut{}$t=1.0$}}}%
      \csname LTb\endcsname
      \put(2651,1800){\rotatebox{16}{\makebox(0,0){\strut{}$t=1.6$}}}%
      \csname LTb\endcsname
      \put(2285,2326){\rotatebox{21}{\makebox(0,0){\strut{}$t=4.0$}}}%
    }%
    \gplbacktext
    \put(0,0){\includegraphics{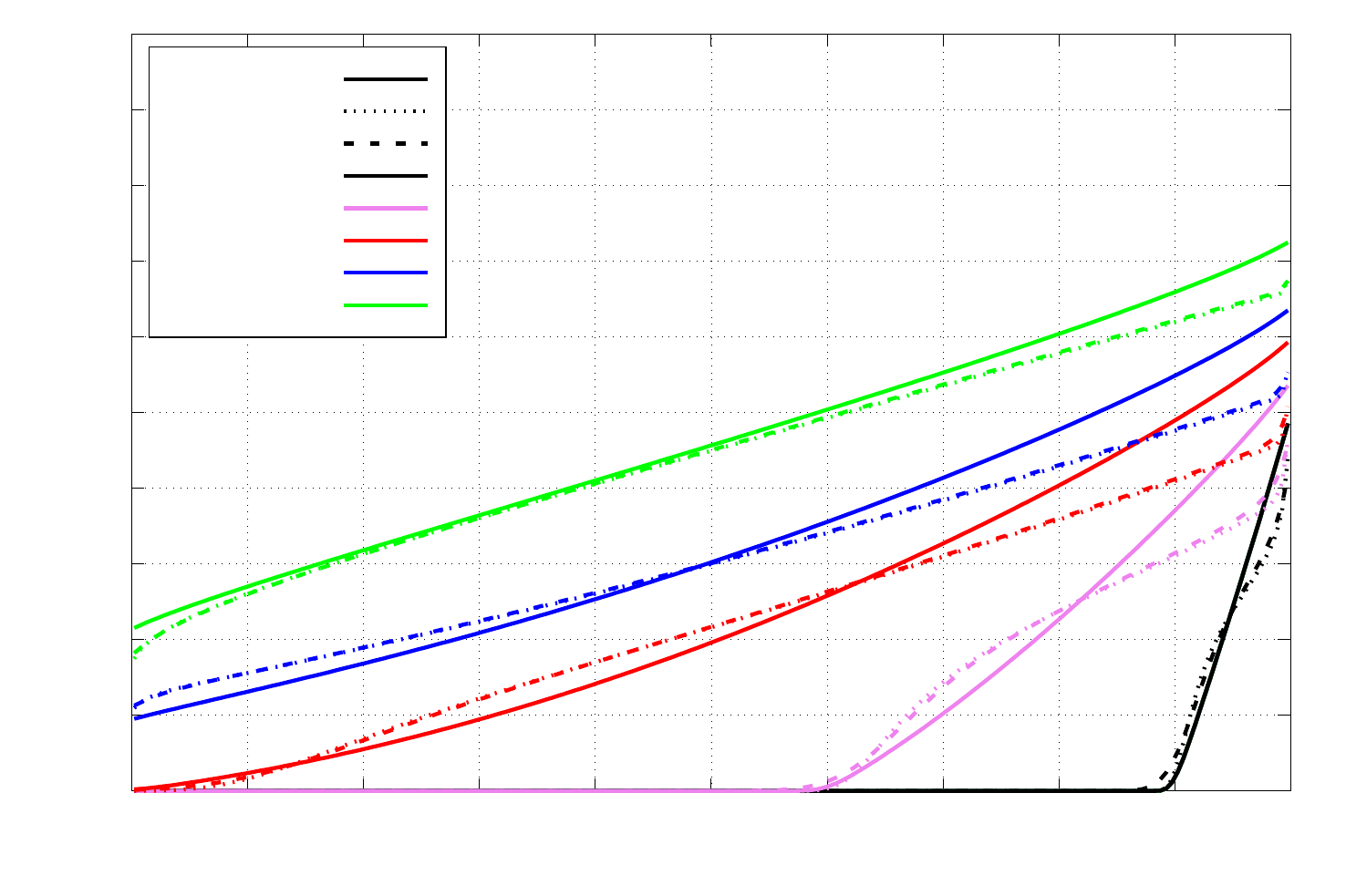}}%
    \gplfronttext
  \end{picture}%
\endgroup

%% file: Aintermediaire.tex
\begingroup
  \makeatletter
  \providecommand\color[2][]{%
    \GenericError{(gnuplot) \space\space\space\@spaces}{%
      Package color not loaded in conjunction with
      terminal option `colourtext'%
    }{See the gnuplot documentation for explanation.%
    }{Either use 'blacktext' in gnuplot or load the package
      color.sty in LaTeX.}%
    \renewcommand\color[2][]{}%
  }%
  \providecommand\includegraphics[2][]{%
    \GenericError{(gnuplot) \space\space\space\@spaces}{%
      Package graphicx or graphics not loaded%
    }{See the gnuplot documentation for explanation.%
    }{The gnuplot epslatex terminal needs graphicx.sty or graphics.sty.}%
    \renewcommand\includegraphics[2][]{}%
  }%
  \providecommand\rotatebox[2]{#2}%
  \@ifundefined{ifGPcolor}{%
    \newif\ifGPcolor
    \GPcolortrue
  }{}%
  \@ifundefined{ifGPblacktext}{%
    \newif\ifGPblacktext
    \GPblacktexttrue
  }{}%
  \let\gplgaddtomacro\g@addto@macro
  \gdef\gplbacktext{}%
  \gdef\gplfronttext{}%
  \makeatother
  \ifGPblacktext
    \def\colorrgb#1{}%
    \def\colorgray#1{}%
  \else
    \ifGPcolor
      \def\colorrgb#1{\color[rgb]{#1}}%
      \def\colorgray#1{\color[gray]{#1}}%
      \expandafter\def\csname LTw\endcsname{\color{white}}%
      \expandafter\def\csname LTb\endcsname{\color{black}}%
      \expandafter\def\csname LTa\endcsname{\color{black}}%
      \expandafter\def\csname LT0\endcsname{\color[rgb]{1,0,0}}%
      \expandafter\def\csname LT1\endcsname{\color[rgb]{0,1,0}}%
      \expandafter\def\csname LT2\endcsname{\color[rgb]{0,0,1}}%
      \expandafter\def\csname LT3\endcsname{\color[rgb]{1,0,1}}%
      \expandafter\def\csname LT4\endcsname{\color[rgb]{0,1,1}}%
      \expandafter\def\csname LT5\endcsname{\color[rgb]{1,1,0}}%
      \expandafter\def\csname LT6\endcsname{\color[rgb]{0,0,0}}%
      \expandafter\def\csname LT7\endcsname{\color[rgb]{1,0.3,0}}%
      \expandafter\def\csname LT8\endcsname{\color[rgb]{0.5,0.5,0.5}}%
    \else
      \def\colorrgb#1{\color{black}}%
      \def\colorgray#1{\color[gray]{#1}}%
      \expandafter\def\csname LTw\endcsname{\color{white}}%
      \expandafter\def\csname LTb\endcsname{\color{black}}%
      \expandafter\def\csname LTa\endcsname{\color{black}}%
      \expandafter\def\csname LT0\endcsname{\color{black}}%
      \expandafter\def\csname LT1\endcsname{\color{black}}%
      \expandafter\def\csname LT2\endcsname{\color{black}}%
      \expandafter\def\csname LT3\endcsname{\color{black}}%
      \expandafter\def\csname LT4\endcsname{\color{black}}%
      \expandafter\def\csname LT5\endcsname{\color{black}}%
      \expandafter\def\csname LT6\endcsname{\color{black}}%
      \expandafter\def\csname LT7\endcsname{\color{black}}%
      \expandafter\def\csname LT8\endcsname{\color{black}}%
    \fi
  \fi
    \setlength{\unitlength}{0.0500bp}%
    \ifx\gptboxheight\undefined%
      \newlength{\gptboxheight}%
      \newlength{\gptboxwidth}%
      \newsavebox{\gptboxtext}%
    \fi%
    \setlength{\fboxrule}{0.5pt}%
    \setlength{\fboxsep}{1pt}%
\begin{picture}(8480.00,5640.00)%
    \gplgaddtomacro\gplbacktext{%
      \csname LTb\endcsname
      \put(708,652){\makebox(0,0)[r]{\strut{}$0$}}%
      \csname LTb\endcsname
      \put(708,1130){\makebox(0,0)[r]{\strut{}$0.1$}}%
      \csname LTb\endcsname
      \put(708,1609){\makebox(0,0)[r]{\strut{}$0.2$}}%
      \csname LTb\endcsname
      \put(708,2087){\makebox(0,0)[r]{\strut{}$0.3$}}%
      \csname LTb\endcsname
      \put(708,2565){\makebox(0,0)[r]{\strut{}$0.4$}}%
      \csname LTb\endcsname
      \put(708,3044){\makebox(0,0)[r]{\strut{}$0.5$}}%
      \csname LTb\endcsname
      \put(708,3522){\makebox(0,0)[r]{\strut{}$0.6$}}%
      \csname LTb\endcsname
      \put(708,4000){\makebox(0,0)[r]{\strut{}$0.7$}}%
      \csname LTb\endcsname
      \put(708,4478){\makebox(0,0)[r]{\strut{}$0.8$}}%
      \csname LTb\endcsname
      \put(708,4957){\makebox(0,0)[r]{\strut{}$0.9$}}%
      \csname LTb\endcsname
      \put(708,5435){\makebox(0,0)[r]{\strut{}$1$}}%
      \csname LTb\endcsname
      \put(820,448){\makebox(0,0){\strut{}$0$}}%
      \csname LTb\endcsname
      \put(1552,448){\makebox(0,0){\strut{}$0.1$}}%
      \csname LTb\endcsname
      \put(2285,448){\makebox(0,0){\strut{}$0.2$}}%
      \csname LTb\endcsname
      \put(3017,448){\makebox(0,0){\strut{}$0.3$}}%
      \csname LTb\endcsname
      \put(3749,448){\makebox(0,0){\strut{}$0.4$}}%
      \csname LTb\endcsname
      \put(4482,448){\makebox(0,0){\strut{}$0.5$}}%
      \csname LTb\endcsname
      \put(5214,448){\makebox(0,0){\strut{}$0.6$}}%
      \csname LTb\endcsname
      \put(5946,448){\makebox(0,0){\strut{}$0.7$}}%
      \csname LTb\endcsname
      \put(6678,448){\makebox(0,0){\strut{}$0.8$}}%
      \csname LTb\endcsname
      \put(7411,448){\makebox(0,0){\strut{}$0.9$}}%
      \csname LTb\endcsname
      \put(8143,448){\makebox(0,0){\strut{}$1$}}%
    }%
    \gplgaddtomacro\gplfronttext{%
      \csname LTb\endcsname
      \put(186,3043){\rotatebox{-270}{\makebox(0,0){\strut{}$\rho(x)$}}}%
      \csname LTb\endcsname
      \put(4481,142){\makebox(0,0){\strut{}$x$}}%
      \csname LTb\endcsname
      \put(2052,5150){\makebox(0,0)[r]{\strut{}UGKS}}%
      \csname LTb\endcsname
      \put(2052,4946){\makebox(0,0)[r]{\strut{}  UGKS-M1}}%
      \csname LTb\endcsname
      \put(2052,4742){\makebox(0,0)[r]{\strut{}HLL}}%
      \csname LTb\endcsname
      \put(2052,4538){\makebox(0,0)[r]{\strut{}  $t=0.1$}}%
      \csname LTb\endcsname
      \put(2052,4334){\makebox(0,0)[r]{\strut{}  $t=0.4$}}%
      \csname LTb\endcsname
      \put(2052,4130){\makebox(0,0)[r]{\strut{}  $t=1.0$}}%
      \csname LTb\endcsname
      \put(2052,3926){\makebox(0,0)[r]{\strut{}  $t=1.6$}}%
      \csname LTb\endcsname
      \put(5946,1585){\rotatebox{42}{\makebox(0,0){\strut{}$t=0.1$}}}%
      \csname LTb\endcsname
      \put(3749,1752){\rotatebox{27}{\makebox(0,0){\strut{}$t=0.4$}}}%
      \csname LTb\endcsname
      \put(3017,1800){\rotatebox{27}{\makebox(0,0){\strut{}$t=1.0$}}}%
      \csname LTb\endcsname
      \put(2285,1943){\rotatebox{31}{\makebox(0,0){\strut{}$t=1.6$}}}%
    }%
    \gplbacktext
    \put(0,0){\includegraphics{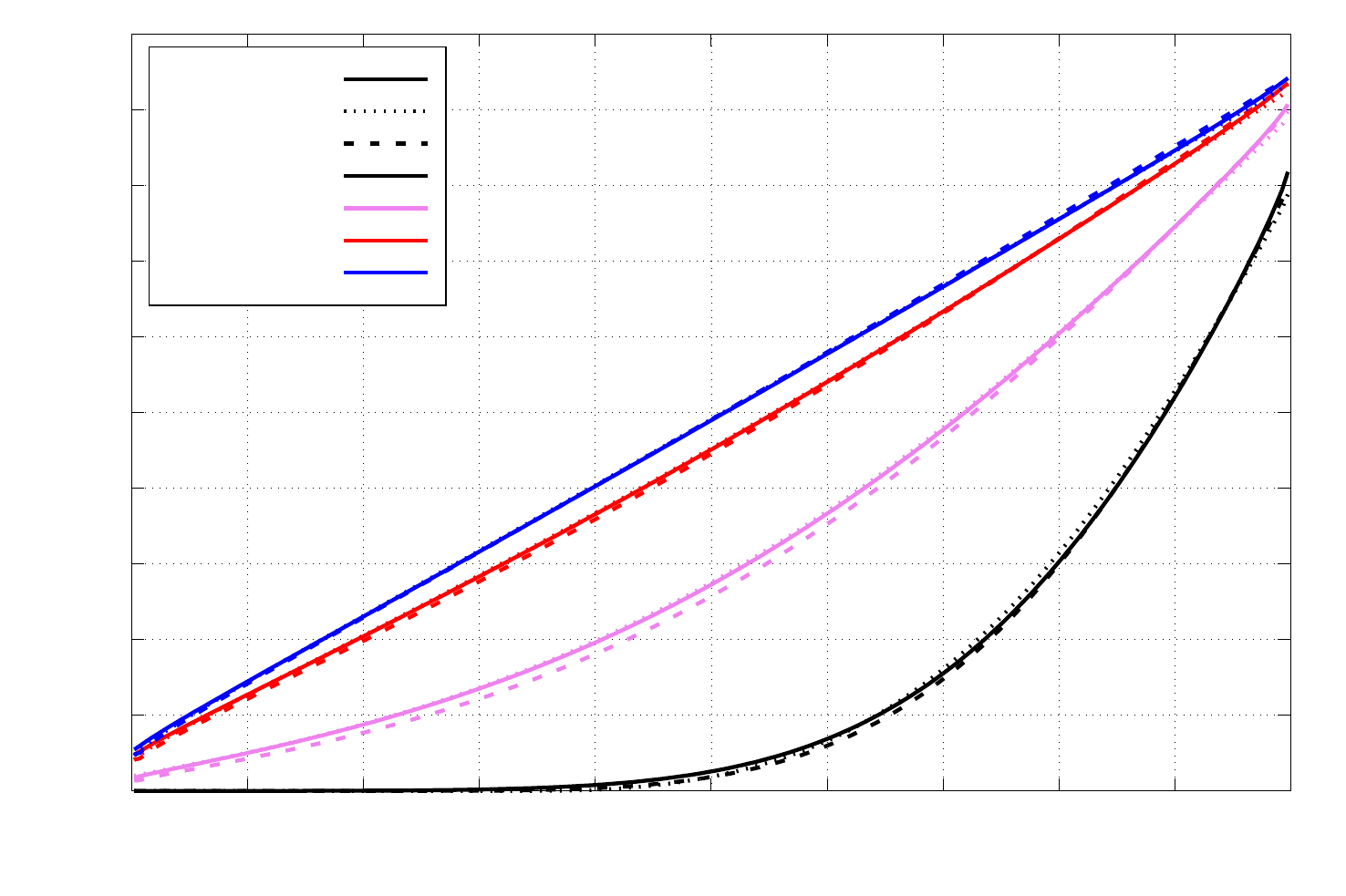}}%
    \gplfronttext
  \end{picture}%
\endgroup

%% file: Adiffusion.tex
\begingroup
  \makeatletter
  \providecommand\color[2][]{%
    \GenericError{(gnuplot) \space\space\space\@spaces}{%
      Package color not loaded in conjunction with
      terminal option `colourtext'%
    }{See the gnuplot documentation for explanation.%
    }{Either use 'blacktext' in gnuplot or load the package
      color.sty in LaTeX.}%
    \renewcommand\color[2][]{}%
  }%
  \providecommand\includegraphics[2][]{%
    \GenericError{(gnuplot) \space\space\space\@spaces}{%
      Package graphicx or graphics not loaded%
    }{See the gnuplot documentation for explanation.%
    }{The gnuplot epslatex terminal needs graphicx.sty or graphics.sty.}%
    \renewcommand\includegraphics[2][]{}%
  }%
  \providecommand\rotatebox[2]{#2}%
  \@ifundefined{ifGPcolor}{%
    \newif\ifGPcolor
    \GPcolortrue
  }{}%
  \@ifundefined{ifGPblacktext}{%
    \newif\ifGPblacktext
    \GPblacktexttrue
  }{}%
  \let\gplgaddtomacro\g@addto@macro
  \gdef\gplbacktext{}%
  \gdef\gplfronttext{}%
  \makeatother
  \ifGPblacktext
    \def\colorrgb#1{}%
    \def\colorgray#1{}%
  \else
    \ifGPcolor
      \def\colorrgb#1{\color[rgb]{#1}}%
      \def\colorgray#1{\color[gray]{#1}}%
      \expandafter\def\csname LTw\endcsname{\color{white}}%
      \expandafter\def\csname LTb\endcsname{\color{black}}%
      \expandafter\def\csname LTa\endcsname{\color{black}}%
      \expandafter\def\csname LT0\endcsname{\color[rgb]{1,0,0}}%
      \expandafter\def\csname LT1\endcsname{\color[rgb]{0,1,0}}%
      \expandafter\def\csname LT2\endcsname{\color[rgb]{0,0,1}}%
      \expandafter\def\csname LT3\endcsname{\color[rgb]{1,0,1}}%
      \expandafter\def\csname LT4\endcsname{\color[rgb]{0,1,1}}%
      \expandafter\def\csname LT5\endcsname{\color[rgb]{1,1,0}}%
      \expandafter\def\csname LT6\endcsname{\color[rgb]{0,0,0}}%
      \expandafter\def\csname LT7\endcsname{\color[rgb]{1,0.3,0}}%
      \expandafter\def\csname LT8\endcsname{\color[rgb]{0.5,0.5,0.5}}%
    \else
      \def\colorrgb#1{\color{black}}%
      \def\colorgray#1{\color[gray]{#1}}%
      \expandafter\def\csname LTw\endcsname{\color{white}}%
      \expandafter\def\csname LTb\endcsname{\color{black}}%
      \expandafter\def\csname LTa\endcsname{\color{black}}%
      \expandafter\def\csname LT0\endcsname{\color{black}}%
      \expandafter\def\csname LT1\endcsname{\color{black}}%
      \expandafter\def\csname LT2\endcsname{\color{black}}%
      \expandafter\def\csname LT3\endcsname{\color{black}}%
      \expandafter\def\csname LT4\endcsname{\color{black}}%
      \expandafter\def\csname LT5\endcsname{\color{black}}%
      \expandafter\def\csname LT6\endcsname{\color{black}}%
      \expandafter\def\csname LT7\endcsname{\color{black}}%
      \expandafter\def\csname LT8\endcsname{\color{black}}%
    \fi
  \fi
    \setlength{\unitlength}{0.0500bp}%
    \ifx\gptboxheight\undefined%
      \newlength{\gptboxheight}%
      \newlength{\gptboxwidth}%
      \newsavebox{\gptboxtext}%
    \fi%
    \setlength{\fboxrule}{0.5pt}%
    \setlength{\fboxsep}{1pt}%
\begin{picture}(8480.00,5640.00)%
    \gplgaddtomacro\gplbacktext{%
      \csname LTb\endcsname
      \put(708,652){\makebox(0,0)[r]{\strut{}$0$}}%
      \csname LTb\endcsname
      \put(708,1130){\makebox(0,0)[r]{\strut{}$0.1$}}%
      \csname LTb\endcsname
      \put(708,1609){\makebox(0,0)[r]{\strut{}$0.2$}}%
      \csname LTb\endcsname
      \put(708,2087){\makebox(0,0)[r]{\strut{}$0.3$}}%
      \csname LTb\endcsname
      \put(708,2565){\makebox(0,0)[r]{\strut{}$0.4$}}%
      \csname LTb\endcsname
      \put(708,3044){\makebox(0,0)[r]{\strut{}$0.5$}}%
      \csname LTb\endcsname
      \put(708,3522){\makebox(0,0)[r]{\strut{}$0.6$}}%
      \csname LTb\endcsname
      \put(708,4000){\makebox(0,0)[r]{\strut{}$0.7$}}%
      \csname LTb\endcsname
      \put(708,4478){\makebox(0,0)[r]{\strut{}$0.8$}}%
      \csname LTb\endcsname
      \put(708,4957){\makebox(0,0)[r]{\strut{}$0.9$}}%
      \csname LTb\endcsname
      \put(708,5435){\makebox(0,0)[r]{\strut{}$1$}}%
      \csname LTb\endcsname
      \put(820,448){\makebox(0,0){\strut{}$0$}}%
      \csname LTb\endcsname
      \put(1552,448){\makebox(0,0){\strut{}$0.1$}}%
      \csname LTb\endcsname
      \put(2285,448){\makebox(0,0){\strut{}$0.2$}}%
      \csname LTb\endcsname
      \put(3017,448){\makebox(0,0){\strut{}$0.3$}}%
      \csname LTb\endcsname
      \put(3749,448){\makebox(0,0){\strut{}$0.4$}}%
      \csname LTb\endcsname
      \put(4482,448){\makebox(0,0){\strut{}$0.5$}}%
      \csname LTb\endcsname
      \put(5214,448){\makebox(0,0){\strut{}$0.6$}}%
      \csname LTb\endcsname
      \put(5946,448){\makebox(0,0){\strut{}$0.7$}}%
      \csname LTb\endcsname
      \put(6678,448){\makebox(0,0){\strut{}$0.8$}}%
      \csname LTb\endcsname
      \put(7411,448){\makebox(0,0){\strut{}$0.9$}}%
      \csname LTb\endcsname
      \put(8143,448){\makebox(0,0){\strut{}$1$}}%
    }%
    \gplgaddtomacro\gplfronttext{%
      \csname LTb\endcsname
      \put(186,3043){\rotatebox{-270}{\makebox(0,0){\strut{}$\rho(x)$}}}%
      \csname LTb\endcsname
      \put(4481,142){\makebox(0,0){\strut{}$x$}}%
      \csname LTb\endcsname
      \put(7278,5150){\makebox(0,0)[r]{\strut{}UGKS}}%
      \csname LTb\endcsname
      \put(7278,4946){\makebox(0,0)[r]{\strut{}  UGKS-M1}}%
      \csname LTb\endcsname
      \put(7278,4742){\makebox(0,0)[r]{\strut{}HLL}}%
      \csname LTb\endcsname
      \put(7278,4538){\makebox(0,0)[r]{\strut{}  $t=0.01$}}%
      \csname LTb\endcsname
      \put(7278,4334){\makebox(0,0)[r]{\strut{}  $t=0.05$}}%
      \csname LTb\endcsname
      \put(7278,4130){\makebox(0,0)[r]{\strut{}  $t=0.15$}}%
      \csname LTb\endcsname
      \put(7278,3926){\makebox(0,0)[r]{\strut{}  $t=2.00$}}%
      \csname LTb\endcsname
      \put(1567,2135){\rotatebox{-72}{\makebox(0,0){\strut{}$t=0.01$}}}%
      \csname LTb\endcsname
      \put(2285,2230){\rotatebox{-60}{\makebox(0,0){\strut{}$t=0.05$}}}%
      \csname LTb\endcsname
      \put(3017,2517){\rotatebox{-46}{\makebox(0,0){\strut{}$t=0.15$}}}%
      \csname LTb\endcsname
      \put(3749,3689){\rotatebox{-33}{\makebox(0,0){\strut{}$t=2.00$}}}%
    }%
    \gplbacktext
    \put(0,0){\includegraphics{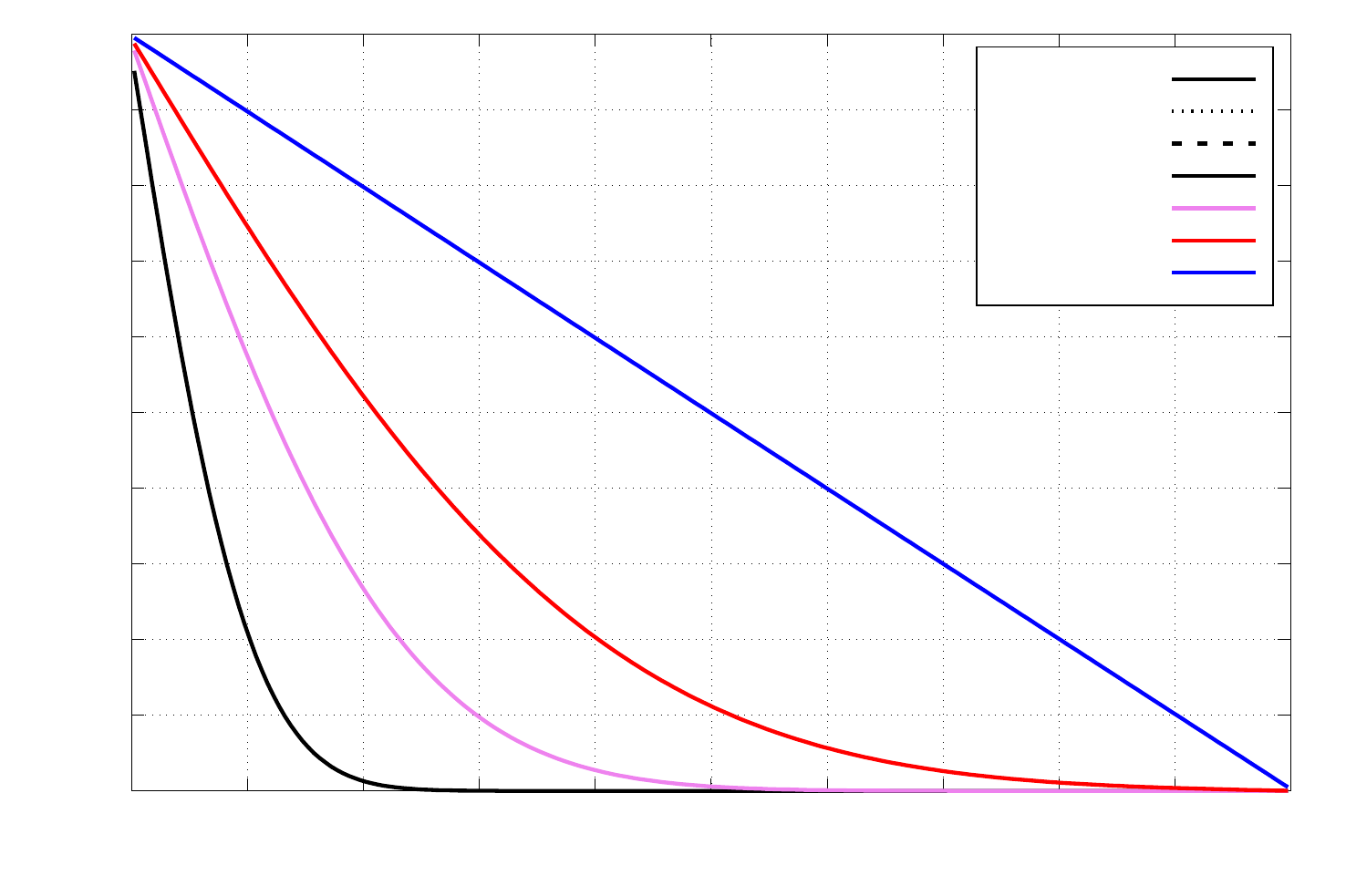}}%
    \gplfronttext
  \end{picture}%
\endgroup

%% file: Acinetique_bis.tex
\begingroup
  \makeatletter
  \providecommand\color[2][]{%
    \GenericError{(gnuplot) \space\space\space\@spaces}{%
      Package color not loaded in conjunction with
      terminal option `colourtext'%
    }{See the gnuplot documentation for explanation.%
    }{Either use 'blacktext' in gnuplot or load the package
      color.sty in LaTeX.}%
    \renewcommand\color[2][]{}%
  }%
  \providecommand\includegraphics[2][]{%
    \GenericError{(gnuplot) \space\space\space\@spaces}{%
      Package graphicx or graphics not loaded%
    }{See the gnuplot documentation for explanation.%
    }{The gnuplot epslatex terminal needs graphicx.sty or graphics.sty.}%
    \renewcommand\includegraphics[2][]{}%
  }%
  \providecommand\rotatebox[2]{#2}%
  \@ifundefined{ifGPcolor}{%
    \newif\ifGPcolor
    \GPcolortrue
  }{}%
  \@ifundefined{ifGPblacktext}{%
    \newif\ifGPblacktext
    \GPblacktexttrue
  }{}%
  \let\gplgaddtomacro\g@addto@macro
  \gdef\gplbacktext{}%
  \gdef\gplfronttext{}%
  \makeatother
  \ifGPblacktext
    \def\colorrgb#1{}%
    \def\colorgray#1{}%
  \else
    \ifGPcolor
      \def\colorrgb#1{\color[rgb]{#1}}%
      \def\colorgray#1{\color[gray]{#1}}%
      \expandafter\def\csname LTw\endcsname{\color{white}}%
      \expandafter\def\csname LTb\endcsname{\color{black}}%
      \expandafter\def\csname LTa\endcsname{\color{black}}%
      \expandafter\def\csname LT0\endcsname{\color[rgb]{1,0,0}}%
      \expandafter\def\csname LT1\endcsname{\color[rgb]{0,1,0}}%
      \expandafter\def\csname LT2\endcsname{\color[rgb]{0,0,1}}%
      \expandafter\def\csname LT3\endcsname{\color[rgb]{1,0,1}}%
      \expandafter\def\csname LT4\endcsname{\color[rgb]{0,1,1}}%
      \expandafter\def\csname LT5\endcsname{\color[rgb]{1,1,0}}%
      \expandafter\def\csname LT6\endcsname{\color[rgb]{0,0,0}}%
      \expandafter\def\csname LT7\endcsname{\color[rgb]{1,0.3,0}}%
      \expandafter\def\csname LT8\endcsname{\color[rgb]{0.5,0.5,0.5}}%
    \else
      \def\colorrgb#1{\color{black}}%
      \def\colorgray#1{\color[gray]{#1}}%
      \expandafter\def\csname LTw\endcsname{\color{white}}%
      \expandafter\def\csname LTb\endcsname{\color{black}}%
      \expandafter\def\csname LTa\endcsname{\color{black}}%
      \expandafter\def\csname LT0\endcsname{\color{black}}%
      \expandafter\def\csname LT1\endcsname{\color{black}}%
      \expandafter\def\csname LT2\endcsname{\color{black}}%
      \expandafter\def\csname LT3\endcsname{\color{black}}%
      \expandafter\def\csname LT4\endcsname{\color{black}}%
      \expandafter\def\csname LT5\endcsname{\color{black}}%
      \expandafter\def\csname LT6\endcsname{\color{black}}%
      \expandafter\def\csname LT7\endcsname{\color{black}}%
      \expandafter\def\csname LT8\endcsname{\color{black}}%
    \fi
  \fi
    \setlength{\unitlength}{0.0500bp}%
    \ifx\gptboxheight\undefined%
      \newlength{\gptboxheight}%
      \newlength{\gptboxwidth}%
      \newsavebox{\gptboxtext}%
    \fi%
    \setlength{\fboxrule}{0.5pt}%
    \setlength{\fboxsep}{1pt}%
\begin{picture}(8480.00,5640.00)%
    \gplgaddtomacro\gplbacktext{%
      \csname LTb\endcsname
      \put(708,652){\makebox(0,0)[r]{\strut{}$0$}}%
      \csname LTb\endcsname
      \put(708,1130){\makebox(0,0)[r]{\strut{}$0.1$}}%
      \csname LTb\endcsname
      \put(708,1609){\makebox(0,0)[r]{\strut{}$0.2$}}%
      \csname LTb\endcsname
      \put(708,2087){\makebox(0,0)[r]{\strut{}$0.3$}}%
      \csname LTb\endcsname
      \put(708,2565){\makebox(0,0)[r]{\strut{}$0.4$}}%
      \csname LTb\endcsname
      \put(708,3044){\makebox(0,0)[r]{\strut{}$0.5$}}%
      \csname LTb\endcsname
      \put(708,3522){\makebox(0,0)[r]{\strut{}$0.6$}}%
      \csname LTb\endcsname
      \put(708,4000){\makebox(0,0)[r]{\strut{}$0.7$}}%
      \csname LTb\endcsname
      \put(708,4478){\makebox(0,0)[r]{\strut{}$0.8$}}%
      \csname LTb\endcsname
      \put(708,4957){\makebox(0,0)[r]{\strut{}$0.9$}}%
      \csname LTb\endcsname
      \put(708,5435){\makebox(0,0)[r]{\strut{}$1$}}%
      \csname LTb\endcsname
      \put(820,448){\makebox(0,0){\strut{}$0$}}%
      \csname LTb\endcsname
      \put(1552,448){\makebox(0,0){\strut{}$0.1$}}%
      \csname LTb\endcsname
      \put(2285,448){\makebox(0,0){\strut{}$0.2$}}%
      \csname LTb\endcsname
      \put(3017,448){\makebox(0,0){\strut{}$0.3$}}%
      \csname LTb\endcsname
      \put(3749,448){\makebox(0,0){\strut{}$0.4$}}%
      \csname LTb\endcsname
      \put(4482,448){\makebox(0,0){\strut{}$0.5$}}%
      \csname LTb\endcsname
      \put(5214,448){\makebox(0,0){\strut{}$0.6$}}%
      \csname LTb\endcsname
      \put(5946,448){\makebox(0,0){\strut{}$0.7$}}%
      \csname LTb\endcsname
      \put(6678,448){\makebox(0,0){\strut{}$0.8$}}%
      \csname LTb\endcsname
      \put(7411,448){\makebox(0,0){\strut{}$0.9$}}%
      \csname LTb\endcsname
      \put(8143,448){\makebox(0,0){\strut{}$1$}}%
    }%
    \gplgaddtomacro\gplfronttext{%
      \csname LTb\endcsname
      \put(186,3043){\rotatebox{-270}{\makebox(0,0){\strut{}$\rho(x)$}}}%
      \csname LTb\endcsname
      \put(4481,142){\makebox(0,0){\strut{}$x$}}%
      \csname LTb\endcsname
      \put(2052,5150){\makebox(0,0)[r]{\strut{}UGKS}}%
      \csname LTb\endcsname
      \put(2052,4946){\makebox(0,0)[r]{\strut{}  UGKS-M1}}%
      \csname LTb\endcsname
      \put(2052,4742){\makebox(0,0)[r]{\strut{}  UGKS-M2}}%
      \csname LTb\endcsname
      \put(2052,4538){\makebox(0,0)[r]{\strut{}  $t=0.1$}}%
      \csname LTb\endcsname
      \put(2052,4334){\makebox(0,0)[r]{\strut{}  $t=0.4$}}%
      \csname LTb\endcsname
      \put(2052,4130){\makebox(0,0)[r]{\strut{}  $t=1.0$}}%
      \csname LTb\endcsname
      \put(2052,3926){\makebox(0,0)[r]{\strut{}  $t=1.6$}}%
      \csname LTb\endcsname
      \put(2052,3722){\makebox(0,0)[r]{\strut{}  $t=4.0$}}%
      \csname LTb\endcsname
      \put(7557,1609){\rotatebox{65}{\makebox(0,0){\strut{}$t=0.1$}}}%
      \csname LTb\endcsname
      \put(5653,1226){\rotatebox{45}{\makebox(0,0){\strut{}$t=0.4$}}}%
      \csname LTb\endcsname
      \put(3017,1369){\rotatebox{20}{\makebox(0,0){\strut{}$t=1.0$}}}%
      \csname LTb\endcsname
      \put(2651,1800){\rotatebox{16}{\makebox(0,0){\strut{}$t=1.6$}}}%
      \csname LTb\endcsname
      \put(2285,2326){\rotatebox{21}{\makebox(0,0){\strut{}$t=4.0$}}}%
    }%
    \gplbacktext
    \put(0,0){\includegraphics{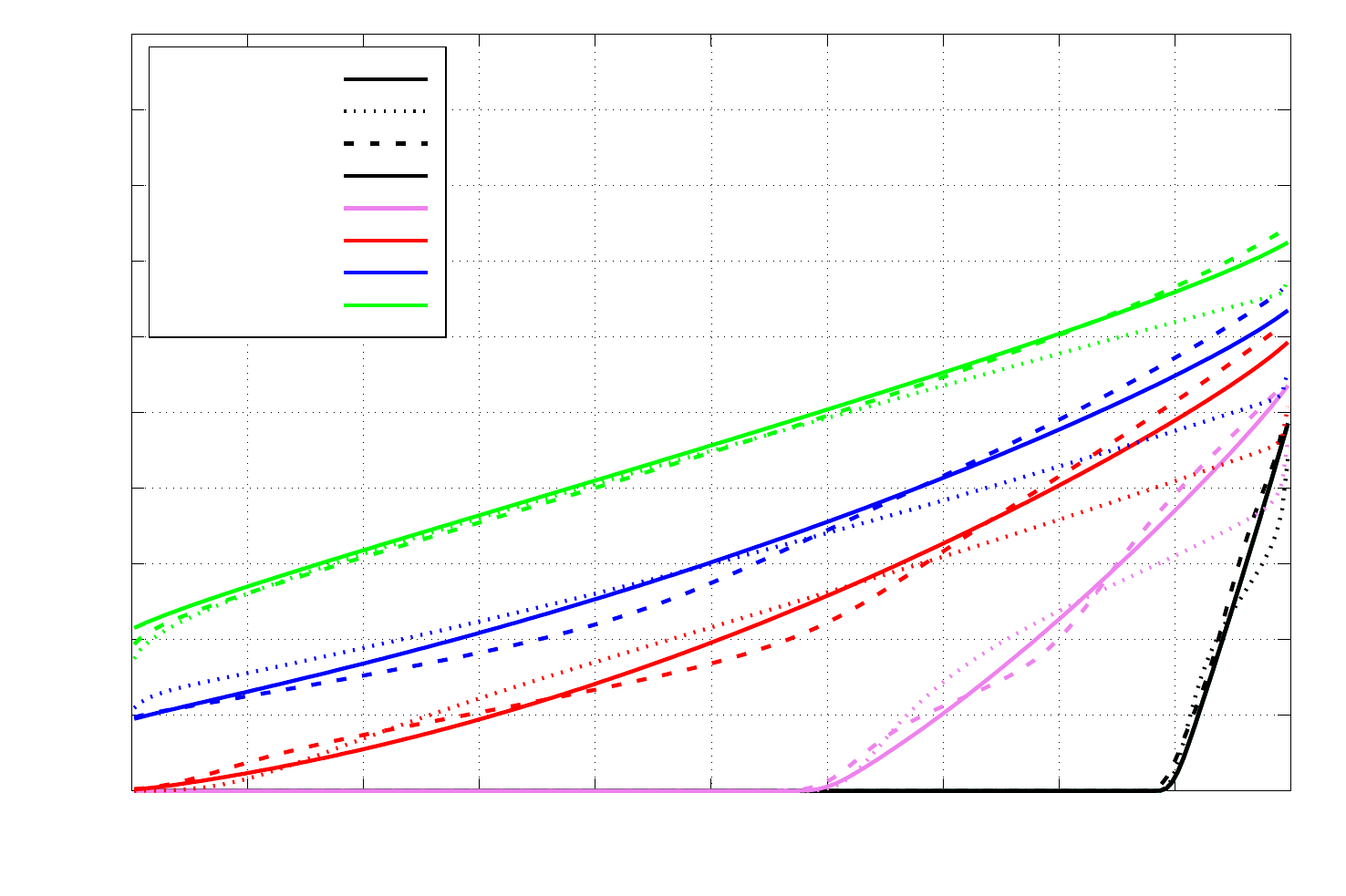}}%
    \gplfronttext
  \end{picture}%
\endgroup

%% file: Aintermediaire_bis.tex
\begingroup
  \makeatletter
  \providecommand\color[2][]{%
    \GenericError{(gnuplot) \space\space\space\@spaces}{%
      Package color not loaded in conjunction with
      terminal option `colourtext'%
    }{See the gnuplot documentation for explanation.%
    }{Either use 'blacktext' in gnuplot or load the package
      color.sty in LaTeX.}%
    \renewcommand\color[2][]{}%
  }%
  \providecommand\includegraphics[2][]{%
    \GenericError{(gnuplot) \space\space\space\@spaces}{%
      Package graphicx or graphics not loaded%
    }{See the gnuplot documentation for explanation.%
    }{The gnuplot epslatex terminal needs graphicx.sty or graphics.sty.}%
    \renewcommand\includegraphics[2][]{}%
  }%
  \providecommand\rotatebox[2]{#2}%
  \@ifundefined{ifGPcolor}{%
    \newif\ifGPcolor
    \GPcolortrue
  }{}%
  \@ifundefined{ifGPblacktext}{%
    \newif\ifGPblacktext
    \GPblacktexttrue
  }{}%
  \let\gplgaddtomacro\g@addto@macro
  \gdef\gplbacktext{}%
  \gdef\gplfronttext{}%
  \makeatother
  \ifGPblacktext
    \def\colorrgb#1{}%
    \def\colorgray#1{}%
  \else
    \ifGPcolor
      \def\colorrgb#1{\color[rgb]{#1}}%
      \def\colorgray#1{\color[gray]{#1}}%
      \expandafter\def\csname LTw\endcsname{\color{white}}%
      \expandafter\def\csname LTb\endcsname{\color{black}}%
      \expandafter\def\csname LTa\endcsname{\color{black}}%
      \expandafter\def\csname LT0\endcsname{\color[rgb]{1,0,0}}%
      \expandafter\def\csname LT1\endcsname{\color[rgb]{0,1,0}}%
      \expandafter\def\csname LT2\endcsname{\color[rgb]{0,0,1}}%
      \expandafter\def\csname LT3\endcsname{\color[rgb]{1,0,1}}%
      \expandafter\def\csname LT4\endcsname{\color[rgb]{0,1,1}}%
      \expandafter\def\csname LT5\endcsname{\color[rgb]{1,1,0}}%
      \expandafter\def\csname LT6\endcsname{\color[rgb]{0,0,0}}%
      \expandafter\def\csname LT7\endcsname{\color[rgb]{1,0.3,0}}%
      \expandafter\def\csname LT8\endcsname{\color[rgb]{0.5,0.5,0.5}}%
    \else
      \def\colorrgb#1{\color{black}}%
      \def\colorgray#1{\color[gray]{#1}}%
      \expandafter\def\csname LTw\endcsname{\color{white}}%
      \expandafter\def\csname LTb\endcsname{\color{black}}%
      \expandafter\def\csname LTa\endcsname{\color{black}}%
      \expandafter\def\csname LT0\endcsname{\color{black}}%
      \expandafter\def\csname LT1\endcsname{\color{black}}%
      \expandafter\def\csname LT2\endcsname{\color{black}}%
      \expandafter\def\csname LT3\endcsname{\color{black}}%
      \expandafter\def\csname LT4\endcsname{\color{black}}%
      \expandafter\def\csname LT5\endcsname{\color{black}}%
      \expandafter\def\csname LT6\endcsname{\color{black}}%
      \expandafter\def\csname LT7\endcsname{\color{black}}%
      \expandafter\def\csname LT8\endcsname{\color{black}}%
    \fi
  \fi
    \setlength{\unitlength}{0.0500bp}%
    \ifx\gptboxheight\undefined%
      \newlength{\gptboxheight}%
      \newlength{\gptboxwidth}%
      \newsavebox{\gptboxtext}%
    \fi%
    \setlength{\fboxrule}{0.5pt}%
    \setlength{\fboxsep}{1pt}%
\begin{picture}(8480.00,5640.00)%
    \gplgaddtomacro\gplbacktext{%
      \csname LTb\endcsname
      \put(708,652){\makebox(0,0)[r]{\strut{}$0$}}%
      \csname LTb\endcsname
      \put(708,1130){\makebox(0,0)[r]{\strut{}$0.1$}}%
      \csname LTb\endcsname
      \put(708,1609){\makebox(0,0)[r]{\strut{}$0.2$}}%
      \csname LTb\endcsname
      \put(708,2087){\makebox(0,0)[r]{\strut{}$0.3$}}%
      \csname LTb\endcsname
      \put(708,2565){\makebox(0,0)[r]{\strut{}$0.4$}}%
      \csname LTb\endcsname
      \put(708,3044){\makebox(0,0)[r]{\strut{}$0.5$}}%
      \csname LTb\endcsname
      \put(708,3522){\makebox(0,0)[r]{\strut{}$0.6$}}%
      \csname LTb\endcsname
      \put(708,4000){\makebox(0,0)[r]{\strut{}$0.7$}}%
      \csname LTb\endcsname
      \put(708,4478){\makebox(0,0)[r]{\strut{}$0.8$}}%
      \csname LTb\endcsname
      \put(708,4957){\makebox(0,0)[r]{\strut{}$0.9$}}%
      \csname LTb\endcsname
      \put(708,5435){\makebox(0,0)[r]{\strut{}$1$}}%
      \csname LTb\endcsname
      \put(820,448){\makebox(0,0){\strut{}$0$}}%
      \csname LTb\endcsname
      \put(1552,448){\makebox(0,0){\strut{}$0.1$}}%
      \csname LTb\endcsname
      \put(2285,448){\makebox(0,0){\strut{}$0.2$}}%
      \csname LTb\endcsname
      \put(3017,448){\makebox(0,0){\strut{}$0.3$}}%
      \csname LTb\endcsname
      \put(3749,448){\makebox(0,0){\strut{}$0.4$}}%
      \csname LTb\endcsname
      \put(4482,448){\makebox(0,0){\strut{}$0.5$}}%
      \csname LTb\endcsname
      \put(5214,448){\makebox(0,0){\strut{}$0.6$}}%
      \csname LTb\endcsname
      \put(5946,448){\makebox(0,0){\strut{}$0.7$}}%
      \csname LTb\endcsname
      \put(6678,448){\makebox(0,0){\strut{}$0.8$}}%
      \csname LTb\endcsname
      \put(7411,448){\makebox(0,0){\strut{}$0.9$}}%
      \csname LTb\endcsname
      \put(8143,448){\makebox(0,0){\strut{}$1$}}%
    }%
    \gplgaddtomacro\gplfronttext{%
      \csname LTb\endcsname
      \put(186,3043){\rotatebox{-270}{\makebox(0,0){\strut{}$\rho(x)$}}}%
      \csname LTb\endcsname
      \put(4481,142){\makebox(0,0){\strut{}$x$}}%
      \csname LTb\endcsname
      \put(2052,5150){\makebox(0,0)[r]{\strut{}UGKS}}%
      \csname LTb\endcsname
      \put(2052,4946){\makebox(0,0)[r]{\strut{}  UGKS-M1}}%
      \csname LTb\endcsname
      \put(2052,4742){\makebox(0,0)[r]{\strut{}  UGKS-M2}}%
      \csname LTb\endcsname
      \put(2052,4538){\makebox(0,0)[r]{\strut{}  $t=0.1$}}%
      \csname LTb\endcsname
      \put(2052,4334){\makebox(0,0)[r]{\strut{}  $t=0.4$}}%
      \csname LTb\endcsname
      \put(2052,4130){\makebox(0,0)[r]{\strut{}  $t=1.0$}}%
      \csname LTb\endcsname
      \put(2052,3926){\makebox(0,0)[r]{\strut{}  $t=1.6$}}%
      \csname LTb\endcsname
      \put(5946,1585){\rotatebox{42}{\makebox(0,0){\strut{}$t=0.1$}}}%
      \csname LTb\endcsname
      \put(3749,1752){\rotatebox{27}{\makebox(0,0){\strut{}$t=0.4$}}}%
      \csname LTb\endcsname
      \put(3017,1800){\rotatebox{27}{\makebox(0,0){\strut{}$t=1.0$}}}%
      \csname LTb\endcsname
      \put(2285,1943){\rotatebox{31}{\makebox(0,0){\strut{}$t=1.6$}}}%
    }%
    \gplbacktext
    \put(0,0){\includegraphics{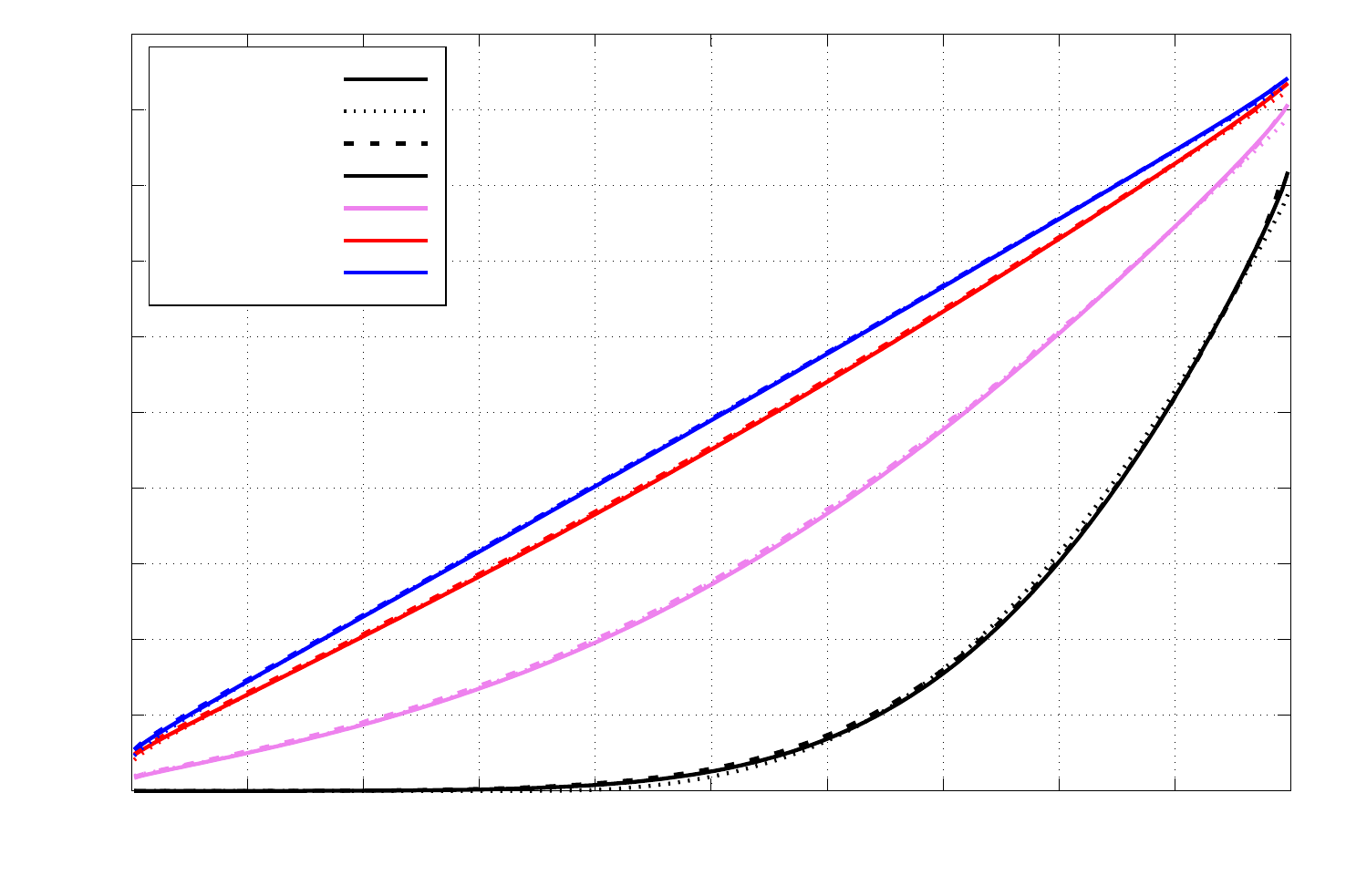}}%
    \gplfronttext
  \end{picture}%
\endgroup

%% file: main.bbl
\begin{thebibliography}{30}
\expandafter\ifx\csname natexlab\endcsname\relax\def\natexlab#1{#1}\fi
\providecommand{\url}[1]{\texttt{#1}}
\providecommand{\href}[2]{#2}
\providecommand{\path}[1]{#1}
\providecommand{\DOIprefix}{doi:}
\providecommand{\ArXivprefix}{arXiv:}
\providecommand{\URLprefix}{URL: }
\providecommand{\Pubmedprefix}{pmid:}
\providecommand{\doi}[1]{\href{http://dx.doi.org/#1}{\path{#1}}}
\providecommand{\Pubmed}[1]{\href{pmid:#1}{\path{#1}}}
\providecommand{\bibinfo}[2]{#2}
\ifx\xfnm\relax \def\xfnm[#1]{\unskip,\space#1}\fi
\bibitem[{Alldredge et~al.(2012)Alldredge, Hauck and Tits}]{alldredge2012high}
\bibinfo{author}{Alldredge, G.W.}, \bibinfo{author}{Hauck, C.D.},
  \bibinfo{author}{Tits, A.L.}, \bibinfo{year}{2012}.
\newblock \bibinfo{title}{High-order entropy-based closures for linear
  transport in slab geometry {II}: A computational study of the optimization
  problem}.
\newblock \bibinfo{journal}{SIAM Journal on Scientific Computing}
  \bibinfo{volume}{34}, \bibinfo{pages}{B361--B391}.
\bibitem[{Bennoune et~al.(2008)Bennoune, Lemou and
  Mieussens}]{bennoune2008uniformly}
\bibinfo{author}{Bennoune, M.}, \bibinfo{author}{Lemou, M.},
  \bibinfo{author}{Mieussens, L.}, \bibinfo{year}{2008}.
\newblock \bibinfo{title}{Uniformly stable numerical schemes for the
  {B}oltzmann equation preserving the compressible {N}avier--{S}tokes
  asymptotics}.
\newblock \bibinfo{journal}{Journal of Computational Physics}
  \bibinfo{volume}{227}, \bibinfo{pages}{3781--3803}.
\bibitem[{Berthon and Turpault(2011)}]{turpault}
\bibinfo{author}{Berthon, C.}, \bibinfo{author}{Turpault, R.},
  \bibinfo{year}{2011}.
\newblock \bibinfo{title}{Asymptotic preserving {HLL} schemes}.
\newblock \bibinfo{journal}{Numerical methods for partial differential
  equations} \bibinfo{volume}{27}, \bibinfo{pages}{1396--1422}.
\bibitem[{Buet et~al.(2002)Buet, Cordier, Lucquin-Desreux and
  Mancini}]{buet2002diffusion}
\bibinfo{author}{Buet, C.}, \bibinfo{author}{Cordier, S.},
  \bibinfo{author}{Lucquin-Desreux, B.}, \bibinfo{author}{Mancini, S.},
  \bibinfo{year}{2002}.
\newblock \bibinfo{title}{Diffusion limit of the {Lorentz} model: asymptotic
  preserving schemes}.
\newblock \bibinfo{journal}{ESAIM: Mathematical Modelling and Numerical
  Analysis} \bibinfo{volume}{36}, \bibinfo{pages}{631--655}.
\bibitem[{Carrillo et~al.(2008a)Carrillo, Goudon and
  Lafitte}]{carrillo2008simulation}
\bibinfo{author}{Carrillo, J.A.}, \bibinfo{author}{Goudon, T.},
  \bibinfo{author}{Lafitte, P.}, \bibinfo{year}{2008}a.
\newblock \bibinfo{title}{Simulation of fluid and particles flows: Asymptotic
  preserving schemes for bubbling and flowing regimes}.
\newblock \bibinfo{journal}{Journal of Computational Physics}
  \bibinfo{volume}{227}, \bibinfo{pages}{7929--7951}.
\bibitem[{Carrillo et~al.(2008b)Carrillo, Goudon, Lafitte and
  Vecil}]{carrillo2008numerical}
\bibinfo{author}{Carrillo, J.A.}, \bibinfo{author}{Goudon, T.},
  \bibinfo{author}{Lafitte, P.}, \bibinfo{author}{Vecil, F.},
  \bibinfo{year}{2008}b.
\newblock \bibinfo{title}{Numerical schemes of diffusion asymptotics and moment
  closures for kinetic equations}.
\newblock \bibinfo{journal}{Journal of Scientific Computing}
  \bibinfo{volume}{36}, \bibinfo{pages}{113--149}.
\bibitem[{Chalons and Guisset(2018)}]{guisset2018}
\bibinfo{author}{Chalons, C.}, \bibinfo{author}{Guisset, S.},
  \bibinfo{year}{2018}.
\newblock \bibinfo{title}{An antidiffusive {HLL} scheme for the electronic
  ${M}_1$ model in the diffusion limit}.
\newblock \bibinfo{journal}{Multiscale Modeling \& Simulation}
  \bibinfo{volume}{16}, \bibinfo{pages}{991--1016}.
\bibitem[{Decoster et~al.(1998)Decoster, Markowich and Perthame}]{modelcoll}
\bibinfo{author}{Decoster, A.}, \bibinfo{author}{Markowich, P.A.},
  \bibinfo{author}{Perthame, B.}, \bibinfo{year}{1998}.
\newblock \bibinfo{title}{Modeling of Collisions}. volume~\bibinfo{volume}{2}.
\newblock \bibinfo{publisher}{Elsevier Masson}.
\bibitem[{Desjardins et~al.(2008)Desjardins, Fox and
  Villedieu}]{desjardins2008quadrature}
\bibinfo{author}{Desjardins, O.}, \bibinfo{author}{Fox, R.O.},
  \bibinfo{author}{Villedieu, P.}, \bibinfo{year}{2008}.
\newblock \bibinfo{title}{A quadrature-based moment method for dilute
  fluid-particle flows}.
\newblock \bibinfo{journal}{Journal of Computational Physics}
  \bibinfo{volume}{227}, \bibinfo{pages}{2514--2539}.
\bibitem[{Dubroca and Feugeas(1999)}]{dubroca1999etude}
\bibinfo{author}{Dubroca, B.}, \bibinfo{author}{Feugeas, J.L.},
  \bibinfo{year}{1999}.
\newblock \bibinfo{title}{Etude th{\'e}orique et num{\'e}rique d'une
  hi{\'e}rarchie de mod{\`e}les aux moments pour le transfert radiatif}.
\newblock \bibinfo{journal}{Comptes Rendus de l'Acad{\'e}mie des
  Sciences-Series I-Mathematics} \bibinfo{volume}{329},
  \bibinfo{pages}{915--920}.
\bibitem[{Gosse(2011)}]{gosse2011transient}
\bibinfo{author}{Gosse, L.}, \bibinfo{year}{2011}.
\newblock \bibinfo{title}{Transient radiative transfer in the grey case:
  Well-balanced and asymptotic-preserving schemes built on case's elementary
  solutions}.
\newblock \bibinfo{journal}{Journal of Quantitative Spectroscopy and Radiative
  Transfer} \bibinfo{volume}{112}, \bibinfo{pages}{1995--2012}.
\bibitem[{Guisset et~al.(2018)Guisset, Brull, Dubroca and
  Turpault}]{guisset2018admissible}
\bibinfo{author}{Guisset, S.}, \bibinfo{author}{Brull, S.},
  \bibinfo{author}{Dubroca, B.}, \bibinfo{author}{Turpault, R.},
  \bibinfo{year}{2018}.
\newblock \bibinfo{title}{An admissible asymptotic-preserving numerical scheme
  for the electronic {M}1 model in the diffusive limit}.
\newblock \bibinfo{journal}{Communications in computational physics}
  \bibinfo{volume}{24}, \bibinfo{pages}{1326--1354}.
\bibitem[{Hauck(2011)}]{hauck2011high}
\bibinfo{author}{Hauck, C.D.}, \bibinfo{year}{2011}.
\newblock \bibinfo{title}{High-order entropy-based closures for linear
  transport in slab geometry}.
\newblock \bibinfo{journal}{Communications in Mathematical Sciences}
  \bibinfo{volume}{9}, \bibinfo{pages}{187--205}.
\bibitem[{Jin and Levermore(1993)}]{jin1993}
\bibinfo{author}{Jin, S.}, \bibinfo{author}{Levermore, C.D.},
  \bibinfo{year}{1993}.
\newblock \bibinfo{title}{Fully-discrete numerical transfer in diffusive
  regimes}.
\newblock \bibinfo{journal}{Transport theory and statistical physics}
  \bibinfo{volume}{22}, \bibinfo{pages}{739--791}.
\bibitem[{Jin and Levermore(1991)}]{jin1991}
\bibinfo{author}{Jin, S.}, \bibinfo{author}{Levermore, D.},
  \bibinfo{year}{1991}.
\newblock \bibinfo{title}{The discrete-ordinate method in diffusive regimes}.
\newblock \bibinfo{journal}{Transport theory and statistical physics}
  \bibinfo{volume}{20}, \bibinfo{pages}{413--439}.
\bibitem[{Jin et~al.(2000)Jin, Pareschi and Toscani}]{jin2000}
\bibinfo{author}{Jin, S.}, \bibinfo{author}{Pareschi, L.},
  \bibinfo{author}{Toscani, G.}, \bibinfo{year}{2000}.
\newblock \bibinfo{title}{Uniformly accurate diffusive relaxation schemes for
  multiscale transport equations}.
\newblock \bibinfo{journal}{SIAM Journal on Numerical Analysis}
  \bibinfo{volume}{38}, \bibinfo{pages}{913--936}.
\bibitem[{Klar(1998)}]{klar1998}
\bibinfo{author}{Klar, A.}, \bibinfo{year}{1998}.
\newblock \bibinfo{title}{An asymptotic-induced scheme for nonstationary
  transport equations in the diffusive limit}.
\newblock \bibinfo{journal}{SIAM journal on numerical analysis}
  \bibinfo{volume}{35}, \bibinfo{pages}{1073--1094}.
\bibitem[{Klar and Schmeiser(2001)}]{klar2001numerical}
\bibinfo{author}{Klar, A.}, \bibinfo{author}{Schmeiser, C.},
  \bibinfo{year}{2001}.
\newblock \bibinfo{title}{Numerical passage from radiative heat transfer to
  nonlinear diffusion models}.
\newblock \bibinfo{journal}{Mathematical Models and Methods in Applied
  Sciences} \bibinfo{volume}{11}, \bibinfo{pages}{749--767}.
\bibitem[{Lafitte and Samaey(2012)}]{lafitte2012asymptotic}
\bibinfo{author}{Lafitte, P.}, \bibinfo{author}{Samaey, G.},
  \bibinfo{year}{2012}.
\newblock \bibinfo{title}{Asymptotic-preserving projective integration schemes
  for kinetic equations in the diffusion limit}.
\newblock \bibinfo{journal}{SIAM Journal on Scientific Computing}
  \bibinfo{volume}{34}, \bibinfo{pages}{A579--A602}.
\bibitem[{Larsen and Morel(1989)}]{larsen1989}
\bibinfo{author}{Larsen, A.W.}, \bibinfo{author}{Morel, J.E.},
  \bibinfo{year}{1989}.
\newblock \bibinfo{title}{Asymptotic solutions of numerical transport problems
  in optically thick, diffusive regimes. {II}}.
\newblock \bibinfo{journal}{Journal of computational physics}
  \bibinfo{volume}{83}, \bibinfo{pages}{212--236}.
\bibitem[{Larsen et~al.(1987)Larsen, Morel and Miller~Jr.}]{larsen1987}
\bibinfo{author}{Larsen, A.W.}, \bibinfo{author}{Morel, J.E.},
  \bibinfo{author}{Miller~Jr., W.F.}, \bibinfo{year}{1987}.
\newblock \bibinfo{title}{Asymptotic solutions of numerical transport problems
  in optically thick, diffusive regimes}.
\newblock \bibinfo{journal}{Journal of computational physics}
  \bibinfo{volume}{69}, \bibinfo{pages}{283--324}.
\bibitem[{Lemou and Mieussens(2008)}]{lemou2008new}
\bibinfo{author}{Lemou, M.}, \bibinfo{author}{Mieussens, L.},
  \bibinfo{year}{2008}.
\newblock \bibinfo{title}{A new asymptotic preserving scheme based on
  micro-macro formulation for linear kinetic equations in the diffusion limit}.
\newblock \bibinfo{journal}{SIAM Journal on Scientific Computing}
  \bibinfo{volume}{31}, \bibinfo{pages}{334--368}.
\bibitem[{Levermore(1996)}]{levermore1996moment}
\bibinfo{author}{Levermore, C.D.}, \bibinfo{year}{1996}.
\newblock \bibinfo{title}{Moment closure hierarchies for kinetic theories}.
\newblock \bibinfo{journal}{Journal of statistical Physics}
  \bibinfo{volume}{83}, \bibinfo{pages}{1021--1065}.
\bibitem[{Liu and Xu(2017)}]{liu2017unified}
\bibinfo{author}{Liu, C.}, \bibinfo{author}{Xu, K.}, \bibinfo{year}{2017}.
\newblock \bibinfo{title}{A unified gas kinetic scheme for continuum and
  rarefied flows {V}: multiscale and multi-component plasma transport}.
\newblock \bibinfo{journal}{Communications in Computational Physics}
  \bibinfo{volume}{22}, \bibinfo{pages}{1175--1223}.
\bibitem[{Mieussens(2013)}]{mieussens2013}
\bibinfo{author}{Mieussens, L.}, \bibinfo{year}{2013}.
\newblock \bibinfo{title}{On the {A}symptotic {P}reserving property of the
  {U}nified {G}as {K}inetic {S}cheme for the diffusion limit of linear kinetic
  models}.
\newblock \bibinfo{journal}{Journal of Computational Physics}
  \bibinfo{volume}{253}, \bibinfo{pages}{138--156}.
\bibitem[{Pichard et~al.(2017)Pichard, Alldredge, Brull, Dubroca and
  Frank}]{pichard2017approximation}
\bibinfo{author}{Pichard, T.}, \bibinfo{author}{Alldredge, G.W.},
  \bibinfo{author}{Brull, S.}, \bibinfo{author}{Dubroca, B.},
  \bibinfo{author}{Frank, M.}, \bibinfo{year}{2017}.
\newblock \bibinfo{title}{An approximation of the m 2 closure: application to
  radiotherapy dose simulation}.
\newblock \bibinfo{journal}{Journal of Scientific Computing}
  \bibinfo{volume}{71}, \bibinfo{pages}{71--108}.
\bibitem[{Sun et~al.(2015)Sun, Jiang, Xu and Li}]{sun2015radiativetransfer}
\bibinfo{author}{Sun, W.}, \bibinfo{author}{Jiang, S.}, \bibinfo{author}{Xu,
  K.}, \bibinfo{author}{Li, S.}, \bibinfo{year}{2015}.
\newblock \bibinfo{title}{An asymptotic preserving unified gas kinetic scheme
  for frequency-dependent radiative transfer equations}.
\newblock \bibinfo{journal}{Journal of Computational Physics}
  \bibinfo{volume}{302}, \bibinfo{pages}{222--238}.
\bibitem[{Van~Leer(1974)}]{vanleer}
\bibinfo{author}{Van~Leer, B.}, \bibinfo{year}{1974}.
\newblock \bibinfo{title}{Towards the ultimate conservative difference scheme.
  {II}. monotonicity and conservation combined in a second-order scheme}.
\newblock \bibinfo{journal}{Journal of computational physics}
  \bibinfo{volume}{14}, \bibinfo{pages}{361--370}.
\bibitem[{Xu and Huang(2010)}]{kxu2010ugks}
\bibinfo{author}{Xu, K.}, \bibinfo{author}{Huang, J.C.}, \bibinfo{year}{2010}.
\newblock \bibinfo{title}{A unified gas-kinetic scheme for continuum and
  rarefied flows}.
\newblock \bibinfo{journal}{Journal of Computational Physics}
  \bibinfo{volume}{229}, \bibinfo{pages}{7747--7764}.
\bibitem[{Zhu and Xu(2021)}]{zhu2021first}
\bibinfo{author}{Zhu, Y.}, \bibinfo{author}{Xu, K.}, \bibinfo{year}{2021}.
\newblock \bibinfo{title}{The first decade of unified gas kinetic scheme}.
\newblock \bibinfo{journal}{arXiv preprint arXiv:2102.01261} .

\end{thebibliography}
